\documentclass{amsart}
\usepackage{amssymb,color,hhline}
\usepackage{amsmath}
\usepackage{a4}
\usepackage{color}
\newtheorem{theorem}{Theorem}[section]
\newtheorem{definition}[theorem]{Definition}
\newtheorem{lemma}[theorem]{Lemma}
\newtheorem{proposition}[theorem]{Proposition}
\newtheorem{observation}[theorem]{Observation}

\newtheorem{remark}[theorem]{Remark}
\newtheorem{corollary}[theorem]{Corollary}

\def\Tr{\operatorname{Tr}}
\def\wedgeo{\wedge}
\def\mapright#1{\smash{\mathop{\longrightarrow}\limits\sp{#1}}}
\makeatletter

 \@addtoreset{equation}{section}
\makeatother
\begin{document}
\title[Geometric theory of equiaffine curvature tensors]
 {Geometric theory of equiaffine curvature tensors}
\author[Gilkey, Nik\v cevi\'c, and Simon]{P. Gilkey, S. Nik\v cevi\'c, and U. Simon }
\begin{address}{PG: Mathematics Department, University of Oregon,
Eugene Or 97403 USA.}
\end{address}
\email{gilkey@uoregon.edu}
\begin{address}
{SN: Mathematical Institute, Sanu, \; 
Knez Mihailova 36, p.p. 367,
 \;  11001 Belgrade, \\
Serbia.}\end{address}
\email{stanan@mi.sanu.ac.rs}
\begin{address}
{US: Institut f\"ur Mathematik, Technische Universit\"at
Berlin \newline  Strasse des 17. Juni 135, D-10623 Berlin, Germany}
\end{address}
\email{simon@math.tu-berlin.de}
%\received{March 30, 2009}
\dedicatory{\begin{center}Dedicated to the memory of Katsumi Nomizu
\end{center}}
\begin{abstract}
From \cite{BGNS} we continue the algebraic investigation of generalized and equiaffine curvature tensors in a given
pseudo-Euclidean vector space  and study different orthogonal, irreducible decompositions in analogy to
the known decomposition of algebraic curvature tensors. We apply
the  decomposition results to characterize geometric properties
of Codazzi structures and relative hypersurfaces; 
particular emphasis is on projectively flat structures.
conjugate connections, projective structures, relative  hypersurfaces.\\
 
\end{abstract}
%%% ----------------------------------------------------------------------
%%% ----------------------------------------------------------------------
%classification:
\maketitle
\section{\bf  Introduction}
%%%%%%%%%%%%%%%%%%%%%%%%%%%%%%%%%%%%%%%%%%%%%%%%%%%%%%%%%%%
In their  famous paper \cite{ST}, I.M. Singer and J.A. Thorpe stated  the
orthogonal decomposition of the Riemannian curvature tensor 
on a 4-manifold into three
components, described by their properties (Definition \ref{defn-6.3}):
\begin{enumerate}
\item constant curvature type,
\item Ricci-traceless,
\item  Ricci-flat. 
\end{enumerate}
This result led to a better understanding of the relations
between algebraic and geometric properties  of the Riemannian 
curvature tensor $R$ and the associated Riemann curvature operator $\mathcal{R}$.
The
studies initiated  a systematic investigation of {\it algebraic curvature tensors}; \cite{G} contains a more complete
bibliography than is possible in this paper. 

Let $V$ be a vector space of dimension $n\ge3$. Let $\mathfrak{a}(V)$
be the space of tensors of type (0,4) with the same symmetries as those of the Riemann curvature tensor. Let
$O(V,g)$ be the orthogonal group associated to a non-degenerate 
scalar  product $g$ on $V$. In Theorem
\ref{thm-6.4}, we will present the well known result that $\mathfrak{a}(V)$ has an irreducible $O(V,g)$
decomposition into the three subspaces described above.

It was Katsumi Nomizu \cite{NOM} who initiated the study of
so called {\it generalized curvature tensors} and 
{\it generalized curvature  operators}, see Definition \ref{defn-2.1}  below;
later other authors, e.g. N. Bokan \cite{BO}, extended his investigations.  Our
paper is devoted to this topic and its geometric applications.

We denote the
real 
 vector space  of  {\it 
generalized curvature operators} by $\mathfrak{R}(V)$. These are
the operators with the same symmetries as the curvature operator of a torsion free connection. Bokan  proved that
the representation of the orthogonal group $O(V,g)$ on
$\mathfrak{R}(V)$ can be decomposed as the sum of eight irreducible subspaces; she dealt with the case that $g$ is
positive definite; 
we refer to \cite{BGNS} for the generalization to arbitrary signatures. This
decomposition is  not unique owing to the fact that two of the representations occur with multiplicity 2 (Lemma
\ref{lem-6.1}).

In  relative  hypersurface theory, in the theory of statistical manifolds,
in the study of Codazzi structures, 
and in Weyl geometry  there appear geometric structures relating to
equiaffine connections, pseudo-Riemannian metrics and
their induced conformal classes. In general their curvature tensors
 do not have the symmetries of the Riemann curvature tensor, but the connections involved  are torsion
free and  admit parallel volume forms.

In \cite{BGNS} we extended known results about the decomposition of generalized curvature 
tensors,   and we 
developed an algebraic
theory of so called {\it equiaffine curvature tensors} (Definition \ref{defn-2.5}). In particular we studied
different orthogonal decompositions of generalized and also of equiaffine curvature tensors for  a given
pseudo-Euclidean vector space.

In this paper we apply our foregoing results and 
develop a geometric theory of generalized and of equiaffine curvature tensors. 
We study  their geometric  properties 
for Codazzi structures with conjugate connections and in relative hypersurface theory. 

For a better understanding of the applications to geometry
it was necessary to extend our  algebraic investigations 
from \cite{BGNS} in more detail in Sections \ref{sect-2} through \ref{sect-6} below. 
For  $(V,g)$  given,
we introduce the following notation: 
the space $\mathfrak{Co}(V)$ is the space of $(1,3)$ curvature operators
satisfying only the standard skew symmetry in the first two arguments; so called {\it generalized curvature operators} additionally 
satisfy the first Bianchi identity; this  space  is denoted by 
$\mathfrak{R}(V)$. In $\mathfrak{Co}(V)$ we introduce 
the concept of {\it $g$-conjugate curvature operators} 
$\mathcal{R}$  and $\mathcal{R}^*$.

The space $\mathfrak{co}(V)$ of
generalized (0,4) curvature tensors
 is $g$-associated to the space  $\mathfrak{Co}(V)$, 
and the space $\mathfrak{r}(V)$  of
generalized (0,4) curvature tensors
is $g$-associated to the space 
$\mathfrak{R}(V)$ of (1,3) curvature operators.
Taking traces with respect to $g$, for $R \in  \mathfrak{r}(V)$  
there appear only two essentially different {\it Ricci type tensors}, denoted by  $Ric$ and 
$Ric^*$;  their role is interchanged by conjugation. 
Both  Ricci type tensors  have the same trace (with respect to the 
scalar product  considered). 
%%%%%%%%%%%%%%%%%%%%%%%%%%%%
%%%%%%%%%%%%%%%%%%%%%%%%%%%%%%%%%%%%%%%%%%%%%%%%%5
We study two different
irreducible, orthogonal decompositions of the space  $\mathfrak{r}(V)$ under the action of the orthogonal group,
each decomposition leads to  eight subspaces: 
$$\mathfrak{r}(V)=W_1\oplus...\oplus W_8=A_1\oplus...\oplus A_8\,.$$
The {\it $W$-decomposition}  induces a decomposition of the space of projective curvature operators. We
will use it subsequently to define additional projective invariants on manifolds.
Similarly, the {\it $A$-decomposition}  induces a decomposition of the space of algebraic curvature tensors. We point
out that  the concept of conjugation of generalized curvature tensors
is a suitable instrument for investigations; this   can be seen from the following 
statement:
\begin{enumerate}
\item We have an  orthogonal $W-$decomposition into 3 subspaces  
$$\mathfrak{r}(V)=W_1  \oplus  [\bigoplus_2^5 W_j]    \oplus  
[\bigoplus_6^8  W_j]\,.$$
\item Any element of $W_1$ is of {\it constant curvature type}.
\item Any element of $\oplus_2^5 W_j$ is {\it Ricci traceless} and 
also $Ricci^*$ {\it traceless}.
\item  Any element of $\oplus_6^8  W_j$ is 
{\it Ricci flat} and  also {\it Ricci${}^*$ flat}.
\end{enumerate}
A similar statement is true for the   $A$-decomposition as we shall discuss presently.
Let $\mathfrak{F}(V) \subset 
\mathfrak{R}(V)$ be the set of equiaffine curvature operators (Definition \ref{defn-2.5}). 
In Observation \ref{obs-6.2}, we discuss  an irreducible, orthogonal  decomposition of $\mathfrak{F}(V)$ into seven 
subspaces.

Our geometric investigations in the second part of the paper
mainly concern the {\it equiaffine}  setting. In many applications we
show how the summands in the two different decompositions 
reflect geometric properties; in particular  we find  new projective invariants.
In the final part we indicate relations
to non-linear PDEs of fourth order that appear as Euler-Lagrange equations
of variational problems in equiaffine hypersurface theory; it is very interesting, that some 
critical points of the Euler-Lagrange equations can be characterized by the vanishing
of some of the components in the   decompositions that we study. 
Applications to geometric structures with non-symmetric Ricci tensors,
in particular to Weyl geometries, shall follow in a subsequent paper.

Here is a brief guide to the paper. The first part of the paper is algebraic in nature. In Section \ref{sect-2}, we
introduce the algebraic theory of curvature tensors and operators, and we present geometric motivations (Theorem
\ref{thm-2.8}). We also discuss 
the {\it conjugate tensor, generalized Ricci tensors}, and
{\it generalized scalar curvatures}. In Section
\ref{sect-3}, we touch briefly on the structure of these spaces as $\operatorname{GL}(V)$ modules. In Section
\ref{sect-4},        we introduce the
$W$-decomposition, and in Section
\ref{sect-5} we introduce the
$A$-decomposition of $\mathfrak{r}(V)$ as $O(V,g)$ modules. 
Some  geometric results are stated in these sections
concerning these decompositions,  and the decompositions are related
to the Ricci and Ricci${}^*$ tensors. 
It is of particular importance that, in the space $\mathfrak{r}(V)$, the
Ricci symmetry of $R$ and $R^*$ is equivalent (Lemma 
\ref{lem-4.9}
and Theorem \ref{thm-4.10}), thus this property is purely algebraic;
so far, a proof was only known in the context of Codazzi structures 
on manifolds in terms of analytic tools (Remark  \ref{rmk-7.4}). 

Let
$\mathfrak{a}(V)$ be the space of algebraic curvature tensors (Definition \ref{defn-2.6}). In Section
\ref{sect-6}, these two decompositions are related and compared to the Singer-Thorpe decomposition of
$\mathfrak{a}(V)$.

The second part of the paper is more geometric in flavor. Section \ref{sect-7} deals with conjugate connections on
manifolds.
Section  \ref{sect-8}  examines Codazzi structures on 
manifolds. Section \ref{sect-9} studies projective
and conformal changes of connections. Section \ref{sect-10} treats relative hypersurface theory. The paper
concludes in Section \ref{sect-11} with an examination of the $W$-decomposition in the framework of relative
hypersurfaces. 

K. Nomizu 
did not only initiate the study of generalized curvature
tensors, he significantly contributed to the
geometry of conjugate connections and affine hypersurface 
theory. 
Our paper treats these topics. We dedicate our investigations
to the memory of this great geometer   of the 20-th century.

%%%%%%%%%%%%%%%%%%%%%%%%%%%%%%%%%%%%%%%%%%%%%%%%%%%%%%%%%%%
\section{Spaces of curvature tensors and operators}\label{sect-2}
In this section we establish notation and provide geometric motivations.

\subsection{Basic Definitions}\label{sect-2.1} Let $V$ be a real vector space of dimension $n$; to simplify the
discussion, we shall assume that $n\ge3$ henceforth. Let $g$ be a
non-degenerate  
scalar product of signature 
$(p,q)$ on $V$.
\begin{definition}\label{defn-2.1}\rm We say that $\mathcal{R}\in\otimes^2V^*\otimes\operatorname{End}(V)$ is a {\it
generalized curvature operator} if it 
satisfies the following  relations
for all $x,y,z\in V$:
\begin{eqnarray}
&&\mathcal{R}(x,y)z=-\mathcal{R}(y,x)z,  \label{eqn-2.a}\\
&&\mathcal{R}(x,y)z+\mathcal{R}(y,z)x+\mathcal{R}(z,x)y=0\,.\label{eqn-2.b}
\end{eqnarray}
As already stated we  denote the space of all  $\mathcal{R}$ satisfying (\ref{eqn-2.a}) by  $\mathfrak{Co}(V),$ and the space of 
generalized curvature operators, satisfying  (\ref{eqn-2.a}) and (\ref{eqn-2.b}),
by $\mathfrak{R}(V)$. Equation 
(\ref{eqn-2.b}) is called the {\it first Bianchi identity}. We use the
scalar product  to raise and lower indices. For $\mathcal{R} \in
\mathfrak{R}(V)$ we define
a corresponding (0,4)-tensor 
$R\in \mathfrak{r}(V)$ 
by means of the identity:
\begin{equation}\label{eqn-2.c}
R(x,y,z,w)=g(\mathcal{R}(x,y)z,w)\,.
\end{equation}
Such a tensor is called a {\it generalized curvature tensor} and is characterized by the identities:
\begin{eqnarray}
&&R(x,y,z,w)=-R(y,x,z,w),\label{eqn-2.d}\\
&&R(x,y,z,w)+R(y,z,x,w)+R(z,x,y,w)=0\,.\label{eqn-2.e}
\end{eqnarray}
Let $\mathfrak{r}(V)$ be the space of all generalized curvature tensors.
The spaces
$\mathfrak{R}(V)$ and $\mathfrak{r}(V)$ are invariant under the action of the general linear group
$\operatorname{GL}(V)$. The isomorphism sending $\mathcal{R}$ to $R$
depends on the scalar product  $g$ or,
equivalently, upon the identification of $V$ with $V^*$; $\mathfrak{R}(V)$ and $\mathfrak{r}(V)$ are not
isomorphic as $\operatorname{GL}(V)$ modules (Remark \ref{rmk-3.2}).
\end{definition}

\begin{definition}\label{defn-2.2}
\rm Let $\mathcal{R}\in\mathfrak{R}(V)$. There are several {\it generalized
Ricci tensors}:
\begin{equation}\label{eqn-2.f}
\begin{array}{l}
\rho_{14}(\mathcal{R})(x,y):=\operatorname{Tr}\{z\rightarrow\mathcal{R}(z,x)y\},\\
\rho_{24}(\mathcal{R})(x,y):=\operatorname{Tr}\{z\rightarrow\mathcal{R}(x,z)y\},\vphantom{\vrule height 11pt}\\
\rho_{34}(\mathcal{R})(x,y):=\operatorname{Tr}\{z\rightarrow\mathcal{R}(x,y)z\}\,;\vphantom{
\vrule height 11pt}
\end{array}\end{equation}
here \; $\Tr$ \; indicates the associated trace operation.
These maps are equivariant with respect to the natural action of $\operatorname{GL}(V)$; there is
no corresponding $\operatorname{GL}(V)$ equivariant map from $\mathfrak{r}(V)$ to $S^2(V^*)$. It follows from
Equations (\ref{eqn-2.a}) and (\ref{eqn-2.b}) that:
\begin{eqnarray*}
&&\rho_{24}(\mathcal{R})(x,y)=-\rho_{14}(\mathcal{R})(x,y),\quad\text{and}\\
&&\rho_{34}(\mathcal{R})(x,y)=-\rho_{14}(\mathcal{R})(x,y)+\rho_{14}(\mathcal{R})(y,x)\,.
\end{eqnarray*}
In particular we have that 
$$\rho_{34}(\mathcal{R}) = 0\quad
\text{if and only if}\quad\rho_{14}(\mathcal{R})\text{ \;\;is symmetric}\,.$$
We adopt the {\it
Einstein convention} and sum over repeated indices.  If $\{e_i\}$ is a basis for $V$, we expand
$\mathcal{R}(e_i,e_j)e_k=\mathcal{R}_{ijk}{}^le_l$, 
$x=x^ie_i$, and $y=y^ie_i$. We then have
\begin{eqnarray*}
\rho_{14}(x,y)=x^iy^j\mathcal{R}_{kij}{}^k,\quad
\rho_{24}(x,y)=x^iy^j\mathcal{R}_{ikj}{}^k,\quad
\rho_{34}(x,y)=x^iy^j\mathcal{R}_{ijk}{}^k\,.
\end{eqnarray*}\end{definition}

\begin{definition}\label{defn-2.3}
\rm Given a {scalar product}  $g$, let $g_{ij}:=g(e_i,e_j)$
and let $g^{ij}$ be the inverse matrix. We use $g$ to define Ricci tensors associated to a generalized
curvature
$R\in\mathfrak{r}(V)$ by setting:
$$\begin{array}{ll}
\rho_{13}(R)(x,y):=g^{ij}R(e_i,x,e_j,x),\;\; &\rho_{14}(R)(x,y):=g^{ij}R(e_i,x,y,e_j),\\
\rho_{23}(R)(x,y):=g^{ij}R(x,e_i,e_j,y),\;\; &\rho_{24}(R)(x,y):=g^{ij}R(x,e_i,y,e_j),\\
\rho_{34}(R)(x,y):=g^{ij}R(x,y,e_i,e_j)\,.
\end{array}$$
\end{definition} 

\begin{definition}\label{defn-2.4}
\rm There is only one relevant scalar geometric invariant  which we shall call
the {\it generalized scalar curvature}
$$\tau:=g^{jk}\mathcal{R}_{ijk}{}^i=g^{il}g^{jk}R_{ijkl}\,.$$
\end{definition}

\begin{definition}\label{defn-2.5}
\rm We say that $\mathcal{F}\in\mathfrak{F}(V)$ is an {\it equiaffine curvature operator} (this notation is
motivated by Definition \ref{defn-2.7}) if, additionally to Equations (\ref{eqn-2.a}) and (\ref{eqn-2.b}), we
have the Ricci  symmetry:
\begin{equation}\label{eqn-2.g}
\rho_{14}(\mathcal{R})(x,y) = \rho_{14}(\mathcal{R})(y,x)\,.
\end{equation}
Let $\mathfrak{F}(V)\subset \mathfrak{R}(V)$ be the
subspace of all equiaffine curvature operators. Again, we use Equation (\ref{eqn-2.c}) to raise indices to define
$\mathfrak{f}(V,g)\subset\mathfrak{r}(V)$; the {scalar product}
$g$ plays a crucial role. The space $\mathfrak{F}(V)$ is a $\operatorname{GL}(V)$ module and the space
$\mathfrak{f}(V,g)$ is an $\operatorname{O}(V,g)$ module.
\end{definition}

\begin{definition}\label{defn-2.6}
\rm The space $\mathfrak{a}(V)\subset\otimes^4V$ of {\it algebraic curvature tensors} is defined by the
following identities:
\begin{eqnarray}
&&A(x,y,z,w)=\; \; {A}(z,w,x,y),\label{eqn-2.h}\\
&&A(x,y,z,w)=-A(y,x,z,w),\label{eqn-2.i}\\
&&A(x,y,z,w)+A(y,z,x,w)+A(z,x,y,w)=0\label{eqn-2.j}\,.
\end{eqnarray}
This space is invariant under the action of $\operatorname{GL}(V)$. If $A\in\mathfrak{a}(V)$ is an algebraic
curvature tensor, then we may use Equation (\ref{eqn-2.c}) to define a corresponding {\it algebraic curvature
operator}
$\mathcal{A}\in\otimes^2V^*\otimes\operatorname{End}(V)$; let
$\mathfrak{A}(V,g)$ be the space of all algebraic curvature operators; this is an $O(V,g)$ module. It is then
immediate that
\begin{eqnarray*}
&&\mathfrak{a}(V)\subset\mathfrak{f}(V,g)\subset\mathfrak{r}(V) \subset \mathfrak{co}(V)  ,\\
&&\mathfrak{A}(V,g)\subset\mathfrak{F}(V)\subset\mathfrak{R}(V) \subset \mathfrak{Co}(V)   \,.
\end{eqnarray*}
\end{definition}

We shall use capital Roman letters $A$, $F$, $R$ for curvature tensors in $\mathfrak{a}(V)$, $\mathfrak{f}(V,g)$,
and
$\mathfrak{r}(V)$, respectively. We shall use capital caligraphic letters $\mathcal{A}$, $\mathcal{F}$, and
$\mathcal{R}$ for the corresponding curvature operators in $\mathfrak{A}(V,g)$, $\mathfrak{F}(V)$, and
$\mathfrak{R}(V)$, respectively. Despite a tendency in the literature to confuse these objects, it is helpful to
distinguish them notationally since the relevant structure groups and module actions differ.

\subsection{Geometric representability I}\label{sect-2.2}
We now present some representability results which provide geometric motivation for our study.
We first establish notation in the geometric setting:

\begin{definition}\label{defn-2.7}
\rm Let  $\nabla$  be a connection on the tangent bundle $TM$ of a
smooth $n$-dimensional manifold $M$.
\begin{enumerate}
\item If $p\in M$, and if $v,w \in T_pM$, the associated {\it
curvature operator} is given by
$$\mathcal{R}_p^\nabla(v,w):=\nabla_v\nabla_w
-\nabla_w\nabla_v-\nabla_{[v,w]}\,.$$
\item If $\nabla$ is torsion free, we say $\nabla$ is {\it equiaffine} if locally there
exists a
$\nabla$-parallel volume element. This is equivalent  to assuming that $\rho_{14}$ is
symmetric \cite{SCHI}.
\item If $g$ is a Riemannian metric on $M$, let $\nabla (g)$ 
be the associated Levi-Civita connection. This
is an equiaffine connection.
\end{enumerate}\end{definition}

%We refer to  We have the following
%representability results \cite{GN03}, that provide the 
%geometric motivation for studying these
%algebraic structures. 
Let $0$ denote the origin of a finite dimensional vector space $V$; if
$\nabla$ is a connection on
$TV$, we let $\mathcal{R}_0^\nabla$ denote the curvature on $T_0V$. 
We have \cite{GN03}:
%%%%%%%%%%%%%%%%%%%%%%%%%%%%%%%%%%%%%%%%%%%%%%%%%%
\begin{theorem}\label{thm-2.8}
\ \begin{enumerate}
\item If $\nabla$ is a torsion free connection on $M$, then $\mathcal{R}_p^\nabla\in\mathfrak{R}(T_pM)$. Conversely,
given $\mathcal{R}\in\mathfrak{R}(V)$, there exists a torsion free connection $\nabla$ on $TV$ so that
$\mathcal{R}_0^\nabla=\mathcal{R}$.
\item If $\nabla$ is an equiaffine connection on $M$, then $\mathcal{R}_p^\nabla\in\mathfrak{F}(T_pM)$.
Conversely, given $\mathcal{F}\in\mathfrak{F}(V)$, there exists an equiaffine connection $\nabla$ on
$TV$ so that $\mathcal{R}_0^\nabla=\mathcal{F}$.
\item We have
$\mathcal{R}_p^{\nabla (g)}\in\mathfrak{A}(T_pM,g_p)$. Conversely, given $\mathcal{R}\in\mathfrak{A}(V,g_0)$, there
exists a pseudo-Riemannian metric $g$ 
on $TV$ so that $g|_{T_0V}=g_0$ and
$\mathcal{R}_0^{\nabla (g)}=\mathcal{R}$.
\end{enumerate}
\end{theorem}
 We postpone until Section \ref{sect-3.1} additional questions
of geometric realizability which arise naturally from the study of $\mathfrak{R}(V)$ as a
$\operatorname{GL}(V)$ module.

\subsection{Conjugation}\label{sect-2.3} 
We return to the algebraic study in $(V,g)$.
The conjugate of a tensor of type (0,4) is defined purely algebraically; to define the conjugate of an operator
requires a {scalar product}. This is a central notion despite the fact that the conjugate of a generalized curvature tensor
(or operator) need not be a generalized curvature tensor (or operator).

\begin{definition}\label{defn-2.9}
\rm Let $R\in\mathfrak{co}(V)$. We define {\it the conjugate tensor}
$$R^*(x,y,z,w):=-R(x,y,w,z)\,.$$
Given a {scalar product} $g$, let $\mathcal{R}$ be the associated curvature operator. Then $\mathcal{R}^*$ is characterized by
the identity:
$$g(\mathcal{R}(x,y)z,w)+g(z,\mathcal{R}^*(x,y)w)=0\,.$$
For this reason, we use the notation {\it conjugate tensor} and {\it conjugate operator} rather than dual tensor
and dual operator. 
\end{definition}
Clearly $R^{**}=R$. 
We observe that, for $R \in \mathfrak{r}(V)$,\quad  $R^* \in \mathfrak{co}(V)$ need not belong to
$\mathfrak{r}(V)$.  In the presence of Equations (\ref{eqn-2.d}) and (\ref{eqn-2.e}),
 Equation (\ref{eqn-2.h}) is equivalent to the identity $A(x,y,z,w)= -A(x,y,w,z)$ \cite{BGNS}. 
Thus we have
$$\mathfrak{a}(V)=\{R\in\mathfrak{r}(V)\; | \; R(x,y,z,w)=-R(x,y,w,z)\}=\{R\in\mathfrak{r}(V)\;
| \; R=R^*\}\,.$$ Thus these tensors are alternating in the last two arguments. It is also
useful to introduce the space of generalized curvature tensors which are symmetric in the last
two  arguments  by setting:
$$\mathfrak{s}(V)=\{R\in\mathfrak{r}(V)\; | \; R(x,y,z,w)=R(x,y,w,z)\}=\{R\in\mathfrak{r}(V) \; | \; R=-R^*\}\,.$$
\begin{lemma}\label{lem-2.10}
 Let $R\in\mathfrak{r}(V)$; then 
$R\in\mathfrak{a}(V) \oplus \mathfrak{s}(V)$ if and only if
$R^*\in\mathfrak{r}(V)$.
\end{lemma}

\begin{proof} If $R\in\mathfrak{a}(V)$, then $R^*=R$. Similarly, if $R\in\mathfrak{s}(V)$, then $R^*=-R$. Thus
if $R\in\mathfrak{a}(V)\oplus\mathfrak{s}(V)$, one has that
$R^*\in\mathfrak{a}(V)\oplus\mathfrak{s}(V)\subset\mathfrak{r}(V)$. This establishes one implication of
the Lemma. Conversely, suppose
$R\in\mathfrak{r}(V)$ and $R^*\in\mathfrak{r}(V)$. We average over the natural $\mathbb{Z}_2$ action interchanging
the last two arguments to define 
$$
R_a:=\textstyle\frac12(R+R^*)\in\mathfrak{a}(V)\quad\text{and}\quad R_s:=\frac12(R-R^*)\in\mathfrak{s}(V)\,.
$$
This shows that
$R=R_a+R_s\in\mathfrak{a}(V)\oplus\mathfrak{s}(V)$ which establishes the other implication of the Lemma. 
\end{proof}
%%%%%%%%%%%%%%%%%%%%%%%%%%%%%%%%%%%%%
We introduce the notation 
$$
Ric(\mathcal{R}):= \rho_{14}(\mathcal{R})\quad\text{and}\quad
Ric^*(\mathcal{R}):= - \rho_{13}(\mathcal{R})= \rho_{23}(\mathcal{R})\,.
$$
We then have that
$$Ric^*(\mathcal{R})=Ric(\mathcal{R}^*)\,.
$$
%%%%%%%%%%%%%%%%%%%%%%%%%%%%%%%%%%%%%
\section{The structure of $\mathfrak{R}(V)$ and $\mathfrak{F}(V)$ as $\operatorname{GL}(V)$ modules}
\label{sect-3}

We have a decomposition of $V^*\otimes V^*$ into irreducible $GL(V)$ modules of the form
$$V^*\otimes V^*=\Lambda^2(V^*)\oplus S^2(V^*)\,.$$
Let $\mathfrak{P}(V):=\ker(\rho_{14})$. Note that
\begin{eqnarray*}
&&\dim\{\mathfrak{P}(V)\}=\textstyle\frac13m^2(m^2-4),\quad
  \dim\{\Lambda^2(V^*)\}=\frac12m(m-1),\\
&&\dim\{S^2(V^*)\}=\textstyle\frac12m(m+1),\quad\dim\{\mathfrak{R}(V)\}=\textstyle\frac13m^2(m^2-1),\\
&&\dim\{\mathfrak{F}(V)\}=\textstyle\frac{m(m-1)(2m^2+2m-3)}6\,.
\end{eqnarray*}
For $\omega\in\Lambda^2(V^*)$ and $\Theta\in S^2(V^*)$, define:
\begin{eqnarray*}
&&\sigma_1(\omega)(x,y)z:=\textstyle\frac{-1}{1+m}\{2\omega(x,y)z+\omega(x,z)y-\omega(y,z)x\},\\
&&\sigma_2(\Theta)(x,y)z:=\textstyle\frac1{1-m}\{\Theta(x,z)y-\Theta(y,z)x\}\,.
\end{eqnarray*}
We refer to Strichartz \cite{S88} for the proof of the following theorem:

\begin{theorem}\label{thm-3.1}
The map $\rho_{14}$ defines two 
$GL(V)$ equivariant short exact sequences
\begin{eqnarray*}
&&0\rightarrow\mathfrak{P}(V)\rightarrow\mathfrak{R}(V)\mapright{\rho_{14}}
\Lambda^2(V^*)\oplus S^2(V^*)\rightarrow0,\\
&&0\rightarrow\mathfrak{P}(V)\rightarrow\mathfrak{F}(V)\mapright{\rho_{14}}
S^2(V^*)\rightarrow 0.
\end{eqnarray*}
which are equivariantly split by the maps $\sigma_1+\sigma_2$ and 
$\sigma_2$, respectively. 
This gives a $GL(V)$ equivariant
decomposition of
\begin{eqnarray*}
&&\mathfrak{R}(V)=\mathfrak{P}(V)\oplus\Lambda^2(V^*)\oplus S^2(V^*),\\
&&\mathfrak{F}(V)\;=\mathfrak{P}(V)\oplus S^2(V^*)
\end{eqnarray*}
as the direct sum of irreducible
$GL(V)$ modules.
\end{theorem}

\begin{remark}\label{rmk-3.2}
\rm Recall that $\dim\mathfrak{a}(V)=m^2(m^2-1)/12$ and that
$\mathfrak{a}(V)$ is an irreducible $\operatorname{GL}(V)$ module
\cite{S88}.
Suppose that $\mathfrak{r}(V)$ and $\mathfrak{R}(V)$ were
isomorphic as $GL(V)$ modules. We would then have $\mathfrak{r}(V)$ as the direct sum of modules of
dimension $m^2(m^2-4)/3$, $m(m+1)/2$, and $m(m-1)/2$ which is impossible. We conclude therefore that the natural
representations of
$\operatorname{GL}(V)$ on $\mathfrak{R}(V)$ and on $\mathfrak{r}(V)$ are not isomorphic. As our primary focus in
this paper is on the $O(V,g)$ module structure, we shall not continue our analysis further of the
$\operatorname{GL}(V)$ module structure of these spaces and instead refer to \cite{BGNS, S88}.
\end{remark}

\subsection{Geometrical representability II}\label{sect-3.1}
There are 8 additional natural geometric realization questions which arise in this context
and whose realizability may be summarized in the following table:
$$\begin{array}{|c|c|c|r||c|c|c|r|}\noalign{\hrule}
\mathfrak{P}(V)&S^2(V^*)&\Lambda^2(V^*)&&\mathfrak{P}(V)&S^2(V^*)&\Lambda^2(V^*)&\\
\noalign{\hrule}\star&\star&\star&\text{yes}&0&\star&\star&\text{yes}\\
\noalign{\hrule}\star&\star&0&\text{yes}&0&\star&0&\text{yes}\\
\noalign{\hrule}\star&0&\star&\text{yes}&0&0&\star&\text{no}\\
\noalign{\hrule}\star&0&0&\text{yes}&0&0&0&\text{yes}\\
\noalign{\hrule}
\end{array}$$
Thus, for example, if 
$$\mathcal{R}(u,v)w = \tfrac{1}{n-1}[Ric(v,w)u - Ric(u,w)v]$$
 and if $Ric(\mathcal{R})$ is
symmetric, then $\mathcal{R}$ can be geometrically realized by a projectively flat,  Ricci
symmetric,  torsion free connection. But if $\mathcal{R}\ne0$ is projectively flat and if
$Ric(\mathcal{R})$ is antisymmetric, then $\mathcal{R}$ can not be geometrically realized by a
projectively flat,  Ricci antisymmetric,  torsion free connection. We refer to \cite{GNW08} for
further details.

\subsection{Rescaling}\label{sect-3.2}

The spaces $\mathfrak{a}(V)$, $\mathfrak{A}(V)$, $\mathfrak{r}(V)$, $\mathfrak{R}(V)$ are $GL(V)$ modules.
We have fixed  a {scalar product}  $g$ on $V$ to raise and lower indices and thereby identify $\mathfrak{R}(V)$ with
$\mathfrak{r}(V)$,  and $\mathfrak{A}(V,g)$ with $\mathfrak{a}(V)$. In terms of components, this isomorphism
may be described by:
\begin{equation}\label{eqn-3.a}
\mathcal{R}_{hij}{}^k \mapsto  \mathcal{R}_{hijl}:=
\mathcal{R}_{hij}{}^k\, g_{kl}\,.
\end{equation}

We can rescale the {scalar product} setting $g_c:=cg$ for
$c>0$.
%this is also sometimes called a {\it } change. 
Thus the isomorphism of Equation
(\ref{eqn-3.a}) has trivial consequences and both, the $A$-decomposition and  the $W$-decomposition, are
unchanged. Such rescalings, however, play a crucial role in invariance theory. H. Weyl's
classical theory of invariance \cite{We46} shows that all $O(V,g)$ scalar invariants of the
curvature tensor (and of its covariant derivatives) arise by contractions of indices. 
The multiplication of a {scalar product}  $g$
on $V$ by a non-zero factor is called a {\it pseudo-conformal change}; 
studying  its effect induces a natural filtration on this space
which is central in many applications. We refer to \cite{K01} for a detailed application of this theory in the
context of heat trace and heat content asymptotics, for example. We also refer to \cite{GKZD} where this analysis
is used to study the graded (or super) trace of the twisted de Rham complex.

\section{The $W$-Decomposition of $\mathfrak{r}(V)$ as an $O(V,g)$ module}\label{sect-4} 

Before stating the first $O(V,g)$ decomposition results for $\mathfrak{r}(V)$, we recall some standard
notation.
%%%%%%%%%%%%%%%%%%%%%%%%%%%%%%%%%%%%%%%%%
\begin{definition}\label{defn-4.1}
\rm
Let $h$ and $k$ be bilinear forms.
\begin{enumerate}
\item Let $ S^2_0(V^*) \subset S^2(V^*)$ be the space of $g-$traceless
symmetric bilinear forms.
\item Set
$h \cdot k(x,y,z,w):= h(x,y)k(z,w)$.
\item For $r=0,1,2,...$, define:
\begin{eqnarray*}
(h \wedge_r k) (x,y,z,w):&=&
h(x,z)k(y,w) - h(y,z)k(x,w)\\
&-&r[h(x,w)k(y,z) - h(y,w)k(x,z)]\,.
\end{eqnarray*}
We set $\wedge:=\wedge_0$ and note that $\wedge_1$ is the    {\it Kulkarni-Nomizu} product:
\begin{eqnarray*}
&&(h \wedge k)(x,y,z,w)= h(x,z)k(y,w) - h(y,z)k(x,w),\\
&&(h \wedge_1 k)(x,y,z,w)=h(x,z)k(y,w) - h(y,z)k(x,w)\\&&\qquad\qquad-h(x,w)k(y,z) + h(y,w)k(x,z)\,.
\end{eqnarray*}
\item Set $\Lambda h (x,y) := \frac{1}{2} [h(x,y) - h(y,x)]$.
\item Set $S h (x,y) := \frac{1}{2} [h(x,y) + h(y,x)]$.
\item Define mappings $\psi$ and $\mu$ from 
$\otimes^4V^*$ to $\otimes^4V^*$  by setting
\begin{eqnarray*}
&&4\psi (R)(x,y,z,w) :=R(x,y,z,w) + R(y,x,w,z)\\
&&\qquad + 
R(z,w,x,y) + R(w,z,y,x);\\
&& 8 \mu (R)(x,y,z,w) := 3 R(x,y,z,w) + 3 R(x,y,w,z)\\
&&\qquad+R(x,w,z,y) +R (x,z,w,y) + R (w,y,z,x) + 
R(z,y,w,x)\,.
\end{eqnarray*}
If we take $R\in\mathfrak{r}(V)$, then $\psi(R)\in\mathfrak{a}(V)$ and $\mu(R)\in\mathfrak{s}(V)$. Furthermore,
$\psi(\psi(R))=\psi(R)$ and $\mu(\mu(R))=\mu(R),$ so these are idempotents \cite{BO}. 
\end{enumerate}
\end{definition}
%%%%%%%%%%%%%%%%%%%%%%%%%%%%%%%%%%%%%%%%%%%%%
\subsection{Components of the $W$-decomposition}\label{sect-4.1} 
We summarize and  extend results from \cite{BGNS,BO}. 
For fixed data $g$ and $R \in  \mathfrak{r}(V)  $, we simply write $Ric:= 
Ric(R),$ $Ric^*: = Ric(R^*),$
and   $\tau:= \tau(R).$ As our calculations are straight forward
we shall omit proofs in the interests of brevity. We may define the
{\it $W$-components} as follows:
\begin{definition}\label{defn-4.2}\rm
Let $\pi_j:\mathfrak{r}(V)\rightarrow W_j$ be the following natural projections:
\medbreak\quad
$\pi_1(R) := \tfrac {- \tau}{n(n-1)} g \wedgeo  g$,
\medbreak\quad
$\pi_2(R) := \tfrac {1}{(n-1)} [\frac{\tau g}{n} -SRic] \wedgeo  g$,
\medbreak\quad
$\pi_3(R) := \tfrac {-1}{(n+1)}  [ 2 \Lambda Ric \cdot g + \Lambda 
Ric \wedgeo  g]$,
\medbreak\quad$\pi_4(R) :=
 \tfrac {-1}{(n^2-4)}[2 \Lambda Ric^{*} \cdot g 
 + 
\Lambda Ric^{*} \wedge_{n+1} g ]$
\smallbreak\qquad\qquad\quad
$-\tfrac {3}{(n^2 - 4)(n+1)} [2 \Lambda Ric  \cdot g 
+ \Lambda Ric \wedge_{n+1} g ]$,
\medbreak\quad
$\pi_5(R) :=  \tfrac {1}{(n-1)(n-2)}[\tau \cdot g \wedgeo  g
-\tfrac{1}{n} S(Ric + (n-1)Ric^*) \wedge_{n-1} g]$,\medbreak\quad
$\pi_6(R) 
:=\psi(R) + \tfrac{1}{2(n-2)} S(Ric + Ric^*)\wedge_1 g
- \tfrac{\tau}{(n-1)(n-2)} g \wedgeo  g$,
\medbreak\quad
$\pi_7(R):= \mu(R) +\tfrac{1}{2n}S(Ric -   Ric^*)\wedge_{-1} g +\tfrac{1}{2(n+2)}\Lambda(3Ric - Ric^*)\cdot g $
\medbreak\qquad\qquad\quad
$+
\tfrac{1}{4(n+2)}\Lambda(3Ric - Ric^*)\wedge_{-1} g,
$
\medbreak\quad
$ \pi_8(R):= R - \psi(R) - \mu(R) +
\tfrac{1}{2(n-2)}\Lambda(Ric +  Ric^*)\cdot g $
$+\tfrac{1}{4(n-2)}\Lambda(Ric +  Ric^*)\wedge_{3} g$.
\end{definition}
\medbreak\noindent  The main result of this section is the
following:
\begin{theorem}\label{thm-4.3} {\bf[W-Decomposition Theorem]}
There is an $O(V,g)$ equivariant orthogonal decomposition of
$\mathfrak{r}(V)=W_1\oplus...\oplus W_8$
as the direct sum of irreducible $O(V,g)$ modules.
\end{theorem}

We note that the isomorphism induced by $g$ identifies $\mathfrak{R}(V)$ with $\mathfrak{r}(V)$ as $O(V,g)$
modules. 
Consequently, Theorem \ref{thm-4.3} also gives the structure of $\mathfrak{R}(V)$ as an $O(V,g)$ module. 
Let
\begin{eqnarray*}
&&\mathfrak{p}(V):=\{R\in\mathfrak{r}(V):Ric(R)=0\},
\\
&&\mathfrak{t}(V):=\{R\in\mathfrak{r}(V):Ric(R)=Ric^*(R)=0\}\subset\mathfrak{p}(V)\,.
\end{eqnarray*}
We have the following characterization of the subspaces $W_j$:
\begin{lemma}\label{lem-4.4}\ 
\begin{enumerate}
\item $R \in W_1$ if and only if $R=cg\wedge g$ for some $c\in\mathbb{R}$.
\item
$R \in W_2$ if and only if $R\in\mathfrak{p}(V)^\perp$ and  $Ric(R)\in S_0^2(V^*)$.
\item $R \in W_3$ if and only if $R\in\mathfrak{p}(V)^\perp$ and $  Ric(R)\in\Lambda^2(V^*)$.
\item $R \in W_4$ if and only if
$R\in \mathfrak{p}(V) \cap   \mathfrak{t}(V)^\perp$ and $Ric^*(R)\in\Lambda^2(V^*)$.
\item  $R \in W_5$ if and only if  
$R\in  \mathfrak{p}(V) \cap \mathfrak{t}(V)^\perp$ and $Ric^*(R)\in S^2(V^*)$.
\item  $R \in W_6$ if and only if $R\in\mathfrak{a}(V)  \cap \mathfrak{t}(V)$.
\item $R \in W_7$ if and only if $
R \in \mathfrak{s}(V) \cap \mathfrak{t}(V)  $.
\item $R \in W_8$ if and only if $R \in (\mathfrak{s}(V) \oplus \mathfrak{a}(V))^{\perp} \cap
\mathfrak{t}(V)  $.
\end{enumerate}\end{lemma}
One may summarize this information in a tabular form.
Denote the projection of $R$ to $\mathfrak{t}(V)$ by $R_o: =R-\sum_{1\le i\le 5}\pi_i(R)$. Then:
\begin{equation}\label{eqn-4.a}
\begin{array}{l}
R_o: = R + \tfrac{2}{n^2 -4} \Lambda[(n-1)Ric +Ric^* ]\cdot g\\
\qquad+\tfrac{1}{n^2 -1}[(n-1) \Lambda Ric +(n+1)S Ric] \wedgeo g\vphantom{\vrule height 11pt}\\
\qquad+\tfrac {1}{(n^2-4)(n+1)} \Lambda(3 Ric + (n+1) Ric^*) \wedge_{n+1} g\vphantom{\vrule height 11pt}\\
\qquad+  \tfrac{1}{n(n-1)(n-2)} S(Ric +(n-1)Ric^*) \wedge_{n-1} g
- \tfrac{\tau}{(n-1)(n-2)} g\wedgeo  g\,.\vphantom{\vrule height 11pt}
\end{array}\end{equation}
\bigbreak\centerline{\bf Table I -- the $W$-decomposition}\medbreak
\centerline{\hbox{$\begin{array}{|l|l|l|}
\noalign{\hrule\hrule}\hfill Ric\ne0\hfill &Ric=0,Ric^*\ne0&Ric=Ric^*=0\\
\noalign{\hrule}\noalign{\hrule}
W_1\hfill(\tau\ne0)&&W_6=\mathfrak{t}(V)\cap\mathfrak{a}(V)\\
\noalign{\hrule}
W_2\ \hfill(Ric\in S_0)&W_5\ \hfill(Ric^*\in S_0)&W_7=
\mathfrak{t}(V)\cap\mathfrak{s}(V)\\
\noalign{\hrule}W_3\hfill(Ric\in\Lambda)&W_4\ \hfill(Ric^*\in
\Lambda)&W_8=\mathfrak{t}(V)\cap\{\mathfrak{a}(V)\oplus\mathfrak{s}(V)\}^\perp\\
\noalign{\hrule}
\end{array}
$}}
\bigbreak\noindent Of course, in  the table we ignore the element
$ 0 \in \mathfrak{r}(V).$
The first column in Table I contains the three components
where the Ricci tensor is non-zero, the second column contains the 2 components where the Ricci
tensor vanishes but the Ricci${}^*$ tensor is non-zero, and the third column contains the 3 components where both,
the Ricci and the Ricci${}^*$ tensors, vanish; thus the third column gives the decomposition of $\mathfrak{t}(V)$.
The first two entries in the third row contain the 2 components where Ricci and Ricci${}^*$ tensors are symmetric
and traceless, and the first two entries in the fourth row contain the 2 components where the Ricci and Ricci${}^*$
tensors are skew symmetric.

The $O(V,g)$ modules $W_i$ are discussed in
\cite{BGNS,BO}. The representations defined by $W_1$, $W_6$, $W_7$, and $W_8$ appear with multiplicity $1$ in
the natural representation of $O(V,g)$ on
$\mathfrak{r}(V)$. 
Thus these summands are unique. On the other hand, the representations corresponding to $W_2$ and $W_5$ are
isomorphic as are the representations corresponding to $W_3$ and $W_4$. Thus these components in
the decomposition of $\mathfrak{r}(V)$ as an $O(V,g)$ module are not unique. This gives rise to the fact that there
can be different decompositions as we shall see when we discuss the
$A$-decomposition in Section \ref{sect-5}. 

$W_6$ is the space of {\it
Weyl conformal curvature tensors}. 
One then has that
$$\pi_6(R) = \psi(R_o),\quad
\pi_7(R)= \mu(R_o),\quad
\pi_8(R)
= R_o - \psi(R_o) - \mu(R_o)\,.$$
%%%%%%%%%%%%%%%%%%%%%%%%%%%%%%%%%%%%%%%%%%%%%%%%%%%%%%%%%%%
\subsection{Properties of the $W$-decomposition}\label{sect-4.2}
A straightforward calculation shows that the Ricci tensors and the Ricci${}^*$ tensors for these components are
given by:
\begin{lemma}\label{lem-4.5} The Ricci tensors of the $W$-components are given by: \begin{enumerate}
\item $Ric(\pi_1(R)) = \tfrac{\tau}{n} g$.
\item $Ric(\pi_2(R)) = - \tfrac{\tau}{n} g + S Ric$.
\item $Ric(\pi_3(R)) = \Lambda Ric$.
\item $Ric(\pi_j(R)) = 0$  for $j=4,...,8$.
\item $Tr_g(Ric(\pi_j(R)) = 0$  for $j=2,...,8$.
\end{enumerate}\end{lemma}
\begin{lemma}\label{lem-4.6}The Ricci${}^*$ tensors of the $W$-components are given by:\begin{enumerate}
\item  $Ric^*(\pi_1(R)) = \tfrac{\tau}{n}  g$.
\item $Ric^*(\pi_2(R)) = \tfrac{1}{n-1}[ \tfrac{\tau}{n} g -  S Ric]$.
\item $Ric^*(\pi_3(R)) = \tfrac{-3}{n+1}\Lambda Ric$.
\item $Ric^*(\pi_4(R)) = \Lambda (Ric^* + \tfrac{3}{n+1}Ric)$.
\item   $Ric^*(\pi_5(R)) =  \tfrac{- \tau}{n-1}g +
S(\tfrac{1}{n-1}Ric + Ric^*)$. 
\item  $Ric^*(\pi_j(R)) = 0$ for $j=6,7,8$.
\item $Tr_g(Ric^*(\pi_j(R)) = 0$ for $j=2,...,8$.
\end{enumerate}\end{lemma}

\begin{lemma}\label{lem-4.7}
The following vanishing results hold:
 \begin{enumerate}
\item $\pi_1(R) = 0$ if and only if $\tau=0$.
\item $\pi_2(R) = 0$  if and only if
$ SRic = \frac{\tau}{n} g$.  
\item $ \pi_3(R) = 0$   if and only if $Ric$ is
symmetric. 
\item $\pi_4(R)=0$ if and only if $\Lambda(Ric^*+\frac3{n+1}Ric)=0$.
\item $\pi_5(R)=0$ if and only if $S(\tfrac{1}{n-1} Ric + Ric^*)
= \frac{\tau}{n-1}g    $.
\end{enumerate}\end{lemma}

We recall the definition of $R_o$ given in Equation (\ref{eqn-4.a}). If $Ric$ and $Ric^*$ are symmetric then the
$W$-components simplify. 
\begin{lemma}\label{lem-4.8}
If $Ric$ and $Ric{}^*$ are symmetric then:
\begin{enumerate}
\item $\pi_1(R) = \tfrac {- \tau}{n(n-1)} g \wedgeo  g$.
\item $\pi_2(R) = \tfrac {1}{(n-1)} [\frac{\tau g}{n} -Ric] \wedgeo  g$.
\item $\pi_3(R)=0$.
\item  $\pi_4(R)=0$. 
\item $\pi_5(R) =  \tfrac {1}{(n-1)(n-2)}[\tau \cdot g \wedgeo  g
-\tfrac{1}{n} (Ric + (n-1)Ric^*) \wedge_{n-1} g]$.
\item 
$R_o = 
R + \frac{1}{n(n -2)}\, (g\,\wedge_{n-1}Ric \, + Ric^* \wedge_{n-1}\,g) -  
\frac{\tau}{(n-1)(n-2)} g \wedgeo  g$.
\item 
$\pi_6(R) 
=\psi(R) + \frac{1}{2(n-2)}(Ric + Ric^*)\wedge_1 g - \frac{\tau}{(n-1)(n-2)} 
g \wedgeo  g$.
\item
$\pi_7(R)
= \mu(R) +\frac{1}{2n}(Ric -   Ric^*)\wedge_{-1} g$.
\item
$
\pi_8(R)(x,y,z,w) = (R - \psi(R) - \mu(R))(x,y,z,w) 
$\smallbreak
$= \tfrac{1}{8}\, \left(3R(x,y,z,w) 
- R(x,y,w,z)
+ R(x,z,w,y)\right)$\smallbreak\qquad
$+  \tfrac{1}{8}\,\left(-3R(x,w,z,y)
+R(z,y,w,x)
+3 R(y,w,z,x)\right).$
\end{enumerate}\end{lemma}
%%%%%%%%%%%%%%%%%%%%%%%%%%%%%%%%%%%%%%%%%%%%%%%%%%%%%%%%%%%%%%%%%%
%%%\begin{lemma}\label{lem-4.9}
%%%If $R\in\mathfrak{f}(V,g)$ and if $R^*\in\mathfrak{f}(V,g)$, then
%%%\begin{enumerate}
%%%\item $\pi_1(R) = \tfrac {- \tau}{n(n-1)} g \wedgeo  g$,
%%%\item $\pi_2(R) = \tfrac {1}{(n-1)} [\frac{\tau g}{n} -Ric] \wedgeo  g$,
%%%\item $\pi_3(R)=0$.
%%%\item  $\pi_4(R)=0$. 
%%%\item $\pi_5(R) =  \tfrac {1}{(n-1)(n-2)}[\tau \cdot g \wedgeo  g
%%%-\tfrac{1}{n} (Ric + (n-1)Ric^*) \wedge_{n-1} g]$
%%%\item 
%%%$\pi_6(R) 
%%%=\psi(R) + \frac{1}{2(n-2)}(Ric + Ric^*)\wedge_1 gg \wedgeo  g$. 
%%%- \frac{\tau}{(n-1)(n-2)} 
%%%\item
%%%$\pi_7(R)
%%%= \mu(R) +\frac{1}{2n}(Ric -   Ric^*)\wedge_{-1} g$.
%%%\item
%%%$\pi_8(R)=0$.
%%%\end{enumerate}\end{lemma}
%%%%%%%%%%%%%%%%%%%%%%%%%%%%%%%%%%%%%%%%%%%%%%%%%%%%%
\begin{lemma}\label{lem-4.9} 
Let  $R\in\mathfrak{r}(V)$ and  $R^*\in\mathfrak{r}(V)$,  then:
\begin{enumerate}
\item $\pi_8(R)= 0 = \pi_8(R^*).$
\item $Ric$ is symmetric if and only if $Ric^*$ is symmetric.
\end{enumerate}
\end{lemma}
\begin{proof} The proof of (2) is elementary, but technical.
From the assumptions we have $R\in  \mathfrak{a}(V) \oplus \mathfrak{s}(V,g)$
and $R^* \in   \mathfrak{a}(V) \oplus \mathfrak{s}(V)$, thus 
 $R =   \mu(R) + \psi(R)$. We insert the definitions of the mappings
$\mu$ and $\psi$ and get:
\begin{eqnarray*}
&&R(x,y,z,w) = \tfrac{1}{4}\left[R(x,y,z,w)+R(y,x,w,z)+R(z,w,x,y)+
R(w,z,y,x)\right] \\
&&\qquad + \tfrac{1}{8}\left[3R(x,y,z,w) + 3R(x,y,w,z) +
R(x,w,z,y)\right] \\
&&\qquad + \tfrac{1}{8}\left[R(x,z,w,y) + R(w,y,z,x) + R(z,y,w,x) \right]\,.
\end{eqnarray*}
Using the skew symmetry and the Bianchi identity, this implies
\begin{eqnarray*}
&&8R(x,y,z,w) =  2\left[R(x,y,z,w)+R(y,x,w,z)+R(z,w,x,y)+R(w,z,y,x)\right]\\
&&\qquad+ 3R(x,y,z,w) + 3R(x,y,w,z) - R(w,z,x,y) - R(z,x,w,y)\\
&&\qquad+ R(x,z,w,y)
- R(y,z,w,x) - R(z,w,y,x) + R(z,y,w,x)\,.
\end{eqnarray*}
We summarize:
\begin{eqnarray*}
3R(x,y,z,w) - R(x,y,w,z) &= &
3R(z,w,x,y) + 2R(x,z,w,y)\\&-& 3R(z,w,y,x)
-2R(y,z,w,x)\,.
\end{eqnarray*}
For the last term use again the Bianchi identity:
\begin{eqnarray*}
&&3R(x,y,z,w)- R(x,y,w,z) \\&=&
3R(z,w,x,y)+ 2R(x,z,w,y)- 3R(z,w,y,x)\\
&+&2R(z,w,y,x) + 2R(w,y,z,x)\,.
\end{eqnarray*}
Take the trace $\Tr\{w \rightarrow R(w,y)z\}$, that yields
\begin{eqnarray*}
&&3Ric(y,z) + Ric^*(y,z) \\
&=& 3Ric^*(z,y) - 2Ric^*(z,y) + 3Ric(z,y) -2Ric(z,y) + 2Ric(y,z)\,,
\end{eqnarray*}
and thus  the identity
$ Ric(y,z) - Ric(z,y) = Ric^*(z,y) - Ric^*(y,z),$ from  which the desired result follows.
\end{proof}
We get the following corollary which we state as a Theorem
according to its importance (see section \ref{sect-10.1.4} below).
\begin{theorem}\label{thm-4.10} 
Let $R, R^*  \in \mathfrak{r}(V)$. Then:
\begin{enumerate}
\item $\mathcal{R}$ is equiaffine { if and only if} $ \mathcal{R}^*$ is equiaffine.
\item
$ R = \pi_1(R) + (\pi_2(R)  + \pi_5(R))   +  (\pi_6(R)  + \pi_7(R)).$
\end{enumerate} 
\end{theorem}
\begin{lemma}\label{lem-4.11} 
Let $R \in \mathfrak{f}(V,g)$ and $R^*  \in \mathfrak{f}(V,g)$; then 
\begin{enumerate}
\item $\pi_1(R) = \tfrac {- \tau}{n(n-1)} g \wedgeo  g = \pi_1(R^*)$. 
\item $\pi_2(R) = \tfrac {1}{(n-1)} [\frac{\tau g}{n} -Ric] \wedgeo  g$.
\item $\pi_3(R)=0 = \pi_3(R^*)$.
\item $\pi_4(R)=0 = \pi_4(R^*)$.
\item $\pi_5(R) =  \tfrac {1}{(n-1)(n-2)}[\tau \cdot g \wedgeo  g
-\tfrac{1}{n} (Ric + (n-1)Ric^*) \wedge_{n-1} g]$.
\item  $\pi_6(R) 
=\psi(R) + \frac{1}{2(n-2)}(Ric + Ric^*)\wedge_1 g - \frac{\tau}{(n-1)(n-2)} 
g \wedgeo  g  = \pi_6(R^*)$.
\item $\pi_7(R)
= \mu(R) +\frac{1}{2n}(Ric -   Ric^*)\wedge_{-1} g = - \pi_7(R^*)$.
\item  $\pi_8(R)=0 = \pi_8(R^*)$.
\end{enumerate}
\end{lemma} 

\begin{remark}\label{lem-4.12}
Let  $R\in\mathfrak{r}(V)$ and  $R^*\in\mathfrak{r}(V)$, then:
$$\mu(R)= \tfrac{1}{2}(R-R^*)\quad\text{and}\quad\psi(R)= \tfrac{1}{2}(R+R^*)\,.$$
\end{remark} 
%%%%%%%%%%%%%%%%%%%%%%%%%%%%%%%%%%%%%%%%%%%%%%%%%%%%%%%%%%%%%%%%%%%%%%%%
\subsection{The Projective Curvature Tensor and Operator}\label{sect-4.3}
We now turn to  projective questions.

\begin{definition}\label{defn-4.13}
\rm Let $g$ be fixed. The {\it projective curvature tensor} $p(R)$ is the projection
of $R$ on $\mathfrak{p}(V)$; this is the space of generalized curvature tensors with $Ric=0$. Thus
$$p(R)  :=\pi_4(R)\oplus...\oplus \pi_8(R) = R - 
[\pi_1(R) \oplus \pi_2(R)  \oplus \pi_3(R)]\,.$$
The $g$-associated (1,3) operator is called the {\it projective curvature operator} and is denoted by
$\mathcal{P}(\mathcal{R})$  \cite{BO}. Note that this definition yields the projective curvature tensor of a
torsion free connection on a manifold. 
\end{definition}

Define
\begin{equation}\label{eqn-4.b}
B^*:=S[Ric^*+(n-1)Ric]-\tau g\,.
\end{equation}

\begin{lemma}\label{lem-4.14} 
Let $R\in\mathfrak{f}(V,g)$. Then
\begin{enumerate}
\item 
$p(R)=\textstyle R + \frac{1}{n-1} (Ric \wedgeo  g)$.
\item If additionally $R^*\in\mathfrak{r}(V)$, then $p(R)=\pi_5(R)+\pi_6(R)+\pi_7(R)$.
Furthermore, the following conditions are equivalent:
\begin{enumerate}
\item
$\mathcal{P}(\mathcal{R}^*) = \mathcal{P}(\mathcal{R})$.
\item $R^* =  R $.
\item $R$ is   algebraic.
\end{enumerate}
\item  If additionally $R^*\in\mathfrak{r}(V)$ and $p(R^*)=0$, then
$B^*=0$.
\end{enumerate}
\end{lemma}
\section{The $A$-Decomposition of $\mathfrak{r}(V)$ as an $O(V,g)$ module} \label{sect-5}
As already stated, the orthogonal decomposition of  $\mathfrak{r}(V)$ 
into irreducible subspaces is not unique. 
In this section we collect and extend results  
on a decomposition different from the  $W$-decomposition. We call it the 
{\it $A$-decomposition}. The 
$A$-components are defined analogously to the $W$-components. In analogy
to Definition \ref{defn-4.2}, we set:
\begin{definition}\label{defn-5.1} Let $\alpha_j:\mathfrak{r}(V)\rightarrow A_i$ be the following
natural projections:
\begin{enumerate}
\item
$\alpha_1({R}) :=\tfrac {- \tau}{n(n-1)}   g \wedgeo  g$. 
\item $\alpha_2({R}) := \tfrac{-1}{2(n-2)} S(Ric +Ric^*)\wedge_1 g
+ \tfrac{2 \tau}{n(n-2)}   g \wedgeo   g$. 
\item $\alpha_3({R}) := \tfrac{-1}{2n} S(Ric - Ric^*)\wedge_{-1} g$.
\item $\alpha_4({R}) := \tfrac{-1}{4(n+2)} [2\Lambda (3 Ric -  Ric^*)g + 
\Lambda (3 Ric -  Ric^*) \wedge_{-1} g]$.
\item $\alpha_5({R}) := \tfrac{-1}{4(n-2)} [2\Lambda (Ric + Ric^*) g + 
\Lambda ( Ric  + Ric^*) \wedge_3 g]$.
\item
$\alpha_6({R}) := \psi({R}) - \alpha_1({R}) 
- \alpha_2({R}) $.
\item
$\alpha_7({R}) := \mu({R}) - \alpha_3({R})
- \alpha_4({R}) $.
\item $\alpha_8({R}) := {R} - \mu({R})
- \psi({R}) - \alpha_5({R})$.
\end{enumerate}
\end{definition}
We have the following analogue of Theorem \ref{thm-4.3}:
\begin{theorem}\label{thm-5.2} 
{\bf[A-Decomposition Theorem]}
There is an $O(V,g)$ equivariant orthogonal
decomposition of
$\mathfrak{r}(V)=A_1\oplus...\oplus A_8$
as the direct sum of orthogonal, irreducible $O(V,g)$ modules.
\end{theorem}
We use Theorem \ref{thm-5.2}
to establish the following useful fact.
\begin{lemma}\label{lem-5.3} 
If $R\in(\mathfrak{a}(V)\oplus\mathfrak{s}(V))^\perp$ and if $Ric(R) \ne 0$, then $Ric$
is skew symmetric and $3Ric=Ric^*$.
\end{lemma}
\begin{proof} Suppose that $R\in(\mathfrak{a}(V)\oplus\mathfrak{s}(V))^\perp$. It then follows that $R$ belongs to
$\ker(\psi)$ and to $\ker(\mu)$. After some calculations
(compare the proof of Lemma \ref{lem-4.9})
this leads to the identity:
\begin{eqnarray*}
0&=&R(x,y,z,w)+2R(x,y,w,z)+R(x,z,w,y)+R(z,y,w,x)\,.
\end{eqnarray*}
Taking trace over $\rho_{13}$ then yields $Ric^*(y,w)=3Ric(y, w)$,  and taking trace over $\rho_{14}$ yields
similarly $Ric(y,z)=-Ric(z,y)$.
\end{proof}
In analogy to Lemma \ref{lem-4.4} we have:
%%%%%%%%%%%%%%%%%%%%%%%%%%%%%%%%%%%%%%%%%%%%%%%%%%
\begin{lemma}\label{lem-5.4}
\ \begin{enumerate}
\item ${R} \in A_1$ if and only if ${R}=cg\wedge g$
for some $c\in\mathbb{R}$.
\item ${R} \in A_2$ if and only if ${R}\in\mathfrak{a}(V) $, $\alpha_6(R)=0$,
and $\tau=0$.
\item ${R} \in A_3$ if and only if ${R} \in 
\mathfrak{s}(V)$, $\alpha_7(R)=0$,
 $Ric({R})\in S^2_0(V^*)$.
\item ${R} \in A_4$ if and only if ${R} \in \mathfrak{s}(V)$, 
$\alpha_7(R)=0$, and $Ric({R})\in\Lambda^2(V^*)$.
\item  ${R} \in A_5$ if and only if ${R} \in (\mathfrak{s}(V) \oplus
\mathfrak{a}(V)\oplus A_8)^{\perp}$.
\item  ${R} \in A_6$ if and only if ${R} \in \mathfrak{a}(V)
\cap \mathfrak{p}(V)$.
\item ${R} \in A_7$ if and only if $ 
{R} \in \mathfrak{s}(V)\cap \mathfrak{p}(V).$  
\item ${R} \in A_8$ if and only if $ 
{R} \in \mathfrak{p}(V) \cap (\mathfrak{s}(V) \oplus \mathfrak{a}(V))^{\perp}$ .
\end{enumerate}
\end{lemma}
%%%%%%%%%%%%%%%%%%%%%%%%%%%%%%%%%%%%%%%%%%%
Again, it is useful to summarize this information in a tabular form:
\medbreak
\centerline{\bf Table II -- the $A$-decomposition}\medbreak
\centerline{\hbox{$\begin{array}{|l||l|l|l|}
\noalign{\hrule}&\mathfrak{a}(V)&\mathfrak{s}(V)&(\mathfrak{a}(V)\oplus\mathfrak{s}(V))^\perp\\
\noalign{\hrule\hrule}\noalign{\hrule}&A_1\hfill(\tau\ne0)&A_4\ \hfill(Ric\in\Lambda)&A_5\ \hfill(Ric\in\Lambda,
3Ric=Ric^*)\\
\noalign{\hrule}&A_2\hfill\ \hfill(Ric\in S_0)&A_3\ \hfill(Ric\in S_0)&\\
\noalign{\hrule}Ric=Ric^*=0&A_6&A_7&A_8\ \hfill\\
\noalign{\hrule}
\end{array}$}}
\bigbreak\noindent 
The first column in Table II contains the components giving the decomposition of
$\mathfrak{a}(V)$, the second column contains the components giving the decomposition of $\mathfrak{s}(V)$, and the
third column contains the components giving the decomposition of $(\mathfrak{a}(V)\oplus\mathfrak{s}(V))^\perp$;
such decompositions are not available from Table I. We also can read off
the symmetry ($S_0$) and skew symmetry ($\Lambda$) of the Ricci tensor from this table. 

\subsection{Properties of the $A$-decomposition}\label{sect-5.1}
Lemmas \ref{lem-4.5} and \ref{lem-4.6}
extend to this setting to become:
\begin{lemma}\label{lem-5.5}
The Ricci tensors of the $A$-components are given by:
\begin{enumerate}
\item $Ric(\alpha_1({R})) = \tfrac{\tau}{n} g$.
\item $Ric(\alpha_2({R})) = \tfrac{-\tau}{n}  g + \tfrac{1}{2} S(Ric + 
Ric^*)$.
\item $Ric(\alpha_3({R})) = \tfrac{1}{2} S(Ric -  Ric^*)$.
\item $Ric(\alpha_4({R})) =  \tfrac{1}{4} \Lambda(3 Ric -  Ric^*)$.
\item $Ric(\alpha_5({R})) = \tfrac{1}{4}  \Lambda(Ric +  Ric^*)$.
\item $ Ric(\alpha_j({R})) = 0$  for  $j = 6,7,8$.
\item $Tr_g(Ric(\alpha_j({R})) = 0$  for  $j=2,...,8$.
\end{enumerate}
\end{lemma}
%%%%%%%%%%%%%%%%%%%%%%%%%%%%%%%%%%%%%%%%%%%%%%%%%
\begin{lemma}\label{lem-5.6}
The Ricci${}^*$ tensors of the $A$-components are given by:
\begin{enumerate}
\item  $Ric^*(\alpha_j({R})) = Ric(\alpha_j({R}))$
 for  $j=1,2$.
\item $Ric^*(\alpha_j({R})) = - Ric(\alpha_j({R}))$
 for  $j=3,4$. 
\item  $Ric^*(\alpha_5({R})) = 3 Ric(\alpha_5({R})).$
\item $Ric^*(\alpha_j({R})) = 0$  for  $j= 6$, $7$, and $8$.
\item $Tr_g(Ric^*(\alpha_j({R})) = 0$  for  $j=2,...,8.$
\end{enumerate}\end{lemma}
%%%%%%%%%%%%%%%%%%%%%%%%%%%%%%%%%%%%%%%%%%%%%%%
As for the $W$-components, there are important vanishing results:
\begin{lemma}\label{lem-5.7}
The following vanishing results hold:\begin{enumerate}
\item $\alpha_1({R}) = 0$  if and only if $\tau=0$.
\item $\alpha_2({R}) = 0$  if and only if
 $ Ric +  Ric^* = \tfrac{2\tau}{n} g$.
\item $ \alpha_3({R}) = 0$   if and only if
 $S Ric =  S Ric^*$.
\item  $ \alpha_4({R}) = 0$ if and only if  $ 3 \Lambda Ric = \Lambda  
Ric^*$.
\item $ \alpha_5({R}) = 0$ if and only if
$\Lambda Ric = - \Lambda  Ric^*$.
\item The following conditions are equivalent:
\begin{enumerate}
\item $R\in\mathfrak{r}(V)$ and $R^*\in\mathfrak{r}(V)$.
\item $R\in\mathfrak{a}(V)\oplus\mathfrak{s}(V)$.
\item $\alpha_5(R)=\alpha_8(R)=0$.
\end{enumerate}\end{enumerate}\end{lemma}
We have as an immediate consequence of Lemma \ref{lem-5.7}(6) that:
\begin{lemma}\label{lem-5.8}
If $R\in\mathfrak{r}(V)$ and if $R^*\in\mathfrak{r}(V)$   then:
\begin{enumerate}
\item $\alpha_j(R)=\alpha_j(R^*)= \alpha_j(R)^*$ for $j=1,2,6$.
\item $ \alpha_3({R^*}) =  (\alpha_3({R}))^* = - \,\alpha_3({R})$.
 \item $ \alpha_j({R^*}) = - \,(\alpha_j({R}))$
 for  $j=4, 7$.  
\item $ \alpha_j({R^*}) = 0 =  \alpha_j({R})$ for $j = 5,8.$   
\end{enumerate}\end{lemma}
%%%%%%%%%%%%%%%%%%%%%%%%%%%%%%%%%%%%%%%
\begin{lemma}\label{lem-5.9}
If $R\in\mathfrak{f}(V,g)$ and if $R^*\in\mathfrak{f}(V,g)$  then:
\begin{enumerate}
\item   $\alpha_1(R)= - \frac{\tau}{n(n-1)} g \wedge g.$
\item    $\alpha_2(R)= - \frac{1}{2(n-2)}(Ric + Ric^*)\wedge_1 g \,+\,\frac{2 \tau}{n(n-2)}\, g \wedge g.$
\item    $\alpha_3(R)= - \frac{1}{2n} (Ric -  Ric^*)\wedge_1 g.$
\item    $\alpha_4(R)= 0 = \alpha_5(R) $.
\item    $\alpha_6(R)= \frac{1}{2}(R + R^*) - \alpha_1(R) - \alpha_2(R).$
\item    $\alpha_7(R)=\frac{1}{2}(R -  R^*) - \alpha_3(R).$
\item    $\alpha_8(R)= 0.$
\end{enumerate}\end{lemma}
%One has the following equivalence result
%\begin{lemma}\label{lem-5.8}
%Let ${R} \in \mathfrak{r}(V)$. 
%The following properties are equivalent:
%\begin{enumerate}
%\item ${R^*} \in \mathfrak{f}(V,g)$.
%\item  $ \alpha_4({R}) = 0$ and  $ \alpha_5({R}) = 0$.
%\item $ \alpha_4({R^*}) = 0$ and  $ \alpha_5({R^*}) = 0$.
%\item $Ric$  and   $Ric^*$ are symmetric.
%\end{enumerate}\end{lemma}
\section{Comparing the $A$-decomposition and the $W$-decomposition}\label{sect-6}
The two decomposition
theorems show that, for given scalar product $g$, there exist different
decompositions of $\mathfrak{r}(V)$
into irreducible and orthogonal subspaces; as noted above, this occurs because not all the representations which
appear have multiplicity one. Referring to the two different decompositions  above, we use the terminology
{\it W-decomposition} and {\it A-decomposition}. 
\begin{lemma}\label{lem-6.1}
\ \begin{enumerate}\item We have the following relations:
$$\begin{array}{lll}
A_1 = W_1,&W_2\oplus W_5=A_2\oplus A_3,&W_3\oplus W_4=A_4\oplus A_5,\\
A_6 = W_6,&A_7 = W_7,&A_8 = W_8\,.\vphantom{\vrule height 12pt}
\end{array}$$
\item As representation spaces of $O(V,g)$, we have isomorphisms
\begin{eqnarray*}
&&W_2\approx W_5\approx A_2\approx A_3,\quad\text{and}\quad
W_3\approx W_4\approx A_4\approx A_5\,.
\end{eqnarray*}
\end{enumerate}\end{lemma}
For the convenience of the reader, we recall once again Tables I and II.
\bigbreak
\centerline{\bf Table I -- the $W$-decomposition}\medbreak
\centerline{\hbox{$\begin{array}{|l|l|l|}
\noalign{\hrule\hrule}\hfill Ric\ne0\hfill &Ric=0,Ric^*\ne0&Ric=Ric^*=0\\
\noalign{\hrule}\noalign{\hrule}
W_1\hfill(\tau\ne0)&&W_6=\mathfrak{t}(V)\cap\mathfrak{a}(V)\\
\noalign{\hrule}
W_2\;\; \hfill(Ric\in S_0)&W_5\ \hfill(Ric^*\in S_0)&W_7=
\mathfrak{t}(V)\cap\mathfrak{s}(V)\\
\noalign{\hrule}W_3\hfill(Ric\in\Lambda)&W_4\ \hfill(Ric^*\in
\Lambda)&W_8=\mathfrak{t}(V)\cap\{\mathfrak{a}(V)\oplus\mathfrak{s}(V)\}^\perp\\
\noalign{\hrule}
\end{array}
$}}
\medbreak
\centerline{\bf Table II -- the $A$-decomposition}\medbreak
\centerline{\hbox{$\begin{array}{|l||l|l|l|}
\noalign{\hrule}&\mathfrak{a}(V)&\mathfrak{s}(V)&(\mathfrak{a}(V)\oplus\mathfrak{s}(V))^\perp\\
\noalign{\hrule\hrule}\noalign{\hrule}&A_1\hfill(\tau\ne0)&A_4\
\hfill(Ric\in\Lambda)&A_5\ \hfill(Ric\in\Lambda,\, 
3Ric=Ric^*)\\
\noalign{\hrule}&A_2\; \hfill(Ric\in S_0)&A_3\ \hfill(Ric\in S_0)&\\
\noalign{\hrule}Ric=Ric^*=0&A_6&A_7&A_8\ \hfill\\
\noalign{\hrule}
\end{array}$}}
\bigbreak\noindent 
Table I contains the decomposition of the space where $Ric=0$; this is the space of projective curvature tensors.
This decomposition is not available from Table II. On the other hand, the components giving the decomposition of
$\mathfrak{a}(V)$, $\mathfrak{s}(V)$, and
$(\mathfrak{a}(V)\oplus\mathfrak{s}(V))^\perp$ are available from Table II but not from Table I.
The elements of the
fourth row in Table II are the same as the elements of the third column in Table I and give the decomposition of
$\mathfrak{t}(V)$.
We summarize this information as follows:
\begin{observation}\label{obs-6.2}
\ \begin{enumerate}
\item The $W$-decomposition allows to recover the subspace of
projective curvature tensors  as the direct sum 
$ \oplus_{j=4}^8 W_j$.  This is not possible in the $A$-decomposition.
\item The $W$-decomposition allows to recover  separately the equiaffine character
of $R$ and  $R^*$, respectively, that means the symmetry 
of $Ric$ and that of $Ric^*$, respectively; that is not possible in
the $A$-decomposition.
\item The $W$-decomposition permits us to express 
$$ \mathfrak{f}(V,g) = W_1 \oplus  W_2  \oplus \bigoplus_{j=4}^8 W_j\,.$$
The $A$-decomposition permits us to express 
$$ \mathfrak{f}(V,g) = \bigoplus_{j=1}^3 A_j \oplus  \bigoplus_{j=6}^8 A_j\,.$$
There are two different decompositions, yet.
\item
In  the decomposition of algebraic curvature tensors
we use the notation {\it constant curvature type}, {\it Ricci traceless}, and {\it Ricci flat};  additionally we shall use 
the notation 
{\it $Ricci^*$ traceless} and {\it $Ricci^*$ flat}.
\item The $A$-decomposition allows us to express $\mathfrak{a}(V)=A_1\oplus A_2\oplus A_6$. This is not possible in
the $W$-decomposition.
\item
We have
$ \oplus_6^8  A_j = \oplus_6^8  W_j$;
 this direct sum is completely traceless and corresponds 
to the Weyl part in the decomposition
of algebraic curvature tensors. 
\end{enumerate}\end{observation}

\subsection{Algebraic curvature tensors}\label{sect-6.1}
We recall the well known  decomposition of algebraic curvature
tensors. 
To compare this decomposition with our decomposition results
in Theorem \ref{thm-6.4} below,  for the convenience of the reader we
recall  the    notation that is often used in the literature (see
e.g. \cite{B}, p. 46):
\begin{definition}\label{defn-6.3}
\rm\ \begin{enumerate}
\item Let $\mathfrak{u}(V,g)$ 
be the space of  all algebraic curvature tensors
of  {\it constant curvature type}. $R\in\mathfrak{u}(V,g)$ if and only if there exists $c\in\mathbb{R}$ so that
$$R(x,y,z,w)=c\{g(x,w)g(y,z)-g(x,z)g(y,w)\}\,.$$
\item Let $\mathfrak{z}(V,g)$ be the space of
all  algebraic curvature tensors that are {\it Ricci traceless}. $R\in\mathfrak{z}(V,g)$ if and only if there exists a
symmetric trace free bilinear form $\Xi$ so that
\begin{eqnarray*}
R(x,y,z,w)&=&\Xi(x,w)g(y,z)+g(x,w)\Xi(y,z)\\&-&\Xi(x,z)g(y,w)-g(x,z)\Xi(y,w)\\&=&-\Xi\wedge_1 g\,.
\end{eqnarray*}
\item Define the subspace  
$\mathfrak{w}(V,g)$ as space of
all  algebraic curvature tensors that are Ricci flat; they are of the type of the {\it Weyl  conformal 
curvature tensor}. We have $R\in\mathfrak{w}(V,g)$ if and only if $R\in\mathfrak{a}(V)$ and
$Ric=0$.
\end{enumerate}\end{definition}
For the comparison, now also recall our notation from Lemma \ref{lem-4.4}.
The following is the celebrated theorem of Singer and Thorpe \cite{ST}.

\begin{theorem}\label{thm-6.4}
{\bf[Algebraic Curvature Decomposition Theorem]} There is an $O(V,g)$ equivariant orthogonal decomposition of
$$\mathfrak{a}(V)=\mathfrak{u}(V,g)\oplus\mathfrak{z}(V,g)\oplus\mathfrak{w}(V,g)$$ 
 as the direct sum
of irreducible and inequivalent $O(V,g)$ modules.
\end{theorem}

The following result compares the decomposition of 
algebraic curvature tensors into 3 orthogonal 
subspaces with the $W$-decomposition and the $A$-decompositions.
\begin{proposition}\label{prop-6.5}
We have the orthogonal decomposition into 3 subspaces  
$$\mathfrak{r}(V)=W_1  \oplus  [\oplus_2^5 W_j]    \oplus
[\oplus_6^8  W_j]
=A_1 \oplus[\oplus_2^5 A_j]\oplus [\oplus_6^8  A_j]\,.$$
\begin{enumerate}
\item  Any $\quad {R} \in \oplus_2^5 W_j=\oplus_2^5A_j \quad$ is Ricci traceless and Ricci${}^*$ traceless.
\item  Any $\quad {R}  \in \oplus_6^8  W_j=\oplus_6^8A_j \quad$ is 
Ricci flat and  Ricci${}^*$ flat.
\item 
$ \mathfrak{u}(V,g) = W_1=A_1, \quad \mathfrak{z}(V,g)=A_2, \quad $ and $\quad \mathfrak{w}(V,g) = A_6= W_6$.
\end{enumerate}\end{proposition}
%%%%%%%%%%%%%%%%%%%%
%%%%%%%%%%%%%%%%%%%%%%%%%%%%%%%%%%%%%%%%%%
%%%%%%%%%%%%%%%%%%%%

\section{Conjugate Connections  on  Manifolds}\label{sect-7}

It is well known that the decomposition of algebraic curvature tensors
reflects geometric properties. 
The following   sections show  that also our foregoing  decomposition results reflect geometric properties.
Let $M$ be a differentiable manifold of dimension $n \ge 3$. We assume
$M$ to be equiped with a  pseudo-Riemannian metric $g$ of signature $(p,q)$.
Let $\nabla$ be a torsion free connection on $M$; in general, $\nabla$ will be different from 
the Levi-Civita connection  $\nabla(g)$. We denote vector fields by 
$u,v,w,...$  Our considerations have
local character. We refer \cite{SSV} for the proofs. 
As already stated in the Introduction, the structure $(\nabla,g)$ appears in many situations.

We say that a smooth differential form $\omega$ of maximal dimension $n$ is a {\it volume form} if $\omega$
is nowhere vanishing; we say that $\omega$ is {\it $\nabla$-parallel}
if $\nabla\omega=0$. 
We summarize   well known facts
for later applications.

\begin{lemma}\label{lem-7.1}
Let $\nabla$ be torsion free.
The following conditions are equivalent:
\begin{enumerate}
\item Locally, $\nabla$ admits a parallel volume form $\omega$.
\item The Ricci tensor 
$Ric:=Ric(\nabla) := Ric(\mathcal{R}(\nabla))$ is symmetric.
\end{enumerate}
\end{lemma}

Note that the local volume form $\omega$ in question is unique modulo a 
non-zero constant factor.

\begin{definition}\label{defn-7.2}
\rm\ \begin{enumerate}
\item A pair $(\nabla, g)$ is called a {\it Codazzi pair} if it satisfies  Codazzi 
equations, that means the covariant derivative $\nabla g$ is totally symmetric.
\item A  triple $(\nabla , g, \nabla ^*)$ of
a non-degenerate metric and two affine connections 
$\nabla$ and $ \nabla ^*$ is called
{\it conjugate} if it satisfies, for all tangent fields $u,v,w,$ the relation 
$u g(v,w) = g(\nabla_u v,w) +  g(v, \nabla^*_u w)$.
Here we admit that the connections $\nabla$ and $ \nabla ^*$ 
have torsion. 
The relation generalizes the well known Ricci Lemma from
Riemannian geometry.
\end{enumerate}\end{definition}
We have the following:
\begin{theorem}\label{thm-7.3}
Let $\nabla$ be an affine  connection on a pseudo-Riemannian manifold $(M,g)$. Then the following
assertions hold:
\begin{itemize}
\item[(A)] \cite{NOM-SI}, \cite{SiS}. Let  $\nabla^\sharp$  be an affine connection.
Then the  triple $(\nabla , g, \nabla ^\sharp)$ is conjugate 
if and only if the pair $(\nabla,g)$ is a Codazzi pair
and the torsion tensors coincide: $T(\nabla) = T(\nabla^\sharp)$. 
\item[(B)] Assume that $\nabla$ is torsion free. Then:
\begin{enumerate}
\item The (1,2) tensor field $C:= \nabla - \nabla (g)$ is symmetric. 
\item The triple $(\nabla, g, \nabla^*)$ 
with $\nabla ^*: = \nabla (g) -  C$   is conjugate.
\item The connection  $\nabla ^*$ is torsion free if and only if 
the pair $(\nabla, g)$ is a Codazzi pair.
\item The connection $\nabla^*$ is Ricci symmetric if and only if the connection $\nabla$ is Ricci symmetric. 
\item The curvature operator $\mathcal{R}$ satisfies the first
Bianchi identity if and only if  $\mathcal{R}^*$ does. 
Thus
$\mathcal{R} \in \mathfrak{r}(T_pM)$ if and only if 
$\mathcal{R}^* \in \mathfrak{r}(T_pM).$
\item  Let  $\nabla$  be  Ricci symmetric. One has that the associated
 curvature  operators $\mathcal{R}=\mathcal{R}(\nabla)$ and $\mathcal{R}^* =\mathcal{R} 
(\nabla^*)$ in a conjugate  triple
are $g-$conjugate.
\end{enumerate}
\end{itemize}
\end{theorem}
\begin{proof} (5) is the only statement without proof in the
literature,  so far. 
Its proof follows from 4.4.10.d in
\cite{SSV}; see also  Lemma \ref{lem-7.6}
and Observation \ref{obs-7.8}
below. 
\end{proof}
\begin{remark}\label{rmk-7.4}
\rm (i)  Note that (4) in Theorem \ref{thm-7.3} needed an analytic
proof
so far (see section 4.4.7 in \cite{SSV}). The proof of Theorem
\ref{thm-4.10}
above is purely algebraic and pointwise.\\
\rm (ii) The sections 4.4.1.i and 4.4.10.d in \cite{SSV} imply also the following: 
Assume that $\nabla$ is torsion free as above;
then  $\nabla^*$ is torsion free if and only if $C$ is symmetric. On the other hand, $\nabla^*$ is torsion free if and only if 
$(\nabla, g)$ is a Codazzi pair. The symmetry of $C$ implies that
$ \gamma _{jkl}{}^i:= C_{jl}{}^h C_{hk}{}^i - C_{jk}{}^h
C_{hl}{}^i$  is algebraic. Thus we can state: Let $\nabla$ be 
torsion free, satisfy the first Bianchi identity, and let 
$(\nabla, g)$ be  a Codazzi pair. Then  $\nabla^*$
satisfies  the first Bianchi identity.\\
\rm (iii)  Consider the conformal class $\mathfrak{C}$ of pseudo-Riemannian
metrics generated by the metric $g$. While the conjugation of a 
connection $\nabla$  in general depends on the choice of a metric in $\mathfrak{C}$,
the conjugation of a generalized curvature tensor in $T_pM$ for $p \in M$ does 
not.\end{remark}

Below   we summarize some facts we shall need concerning the {\it
  cubic form} $C$ and the  {\it Tchebychev form}  $T$:
\begin{observation}\label{obs-7.5}
\rm Let $(\nabla , g, \nabla ^*)$ be a conjugate triple
with torsion free connections $\nabla , \nabla ^*$.
Then:
\begin{enumerate}
\item We have  $C= \frac {1}{2}(\nabla - \nabla^*) =
\nabla(g) - \nabla^* = \nabla - \nabla(g)$.
The trace $T^{\flat}$ of $C$ is a 1-form, called the 
{\it Tchebychev form}, it is given by
$$nT^{\flat}(v):=\Tr (u \mapsto C(u,v)).$$
The $g$-associated  {\it Tchebychev vector field}  $T$ is
characterized by the identity  $g(T,v):= T^{\flat}(v)$; as usual, we use a  simplified 
coordinate notation for $T^{\flat}$ and write $T_i:=T^{\flat}_i$.
\item If the connections $\nabla$ and $\nabla^*$
are Ricci symmetric then they locally 
admit parallel volume forms $\omega$ and $ \omega^*$, respectively; then 
the
Tchebychev form is closed and satisfies
$$nT^{\flat} = d \ln \frac{\omega(v_1,...,v_n)}{\omega(g)(v_1,...,v_n)} 
= d \ln \frac{\omega(g)(v_1,...,v_n)}{\omega^*(v_1,...,v_n)} ,$$
where the volume forms have the same orientation and are evaluated
on the same frame $(v_1,...,v_n)$.
\item Let $ C^{\flat}(u,v,w) := g(C(u,v),w)$ be the totally symmetric
{\it  cubic form} generated by $C$. Its symmetry is
equivalent to the  Codazzi properties of $g$; namely, its covariant derivatives are given by:
$ \nabla^* g = 2 C^{\flat} = - \nabla g$.
\end{enumerate}\end{observation}

For the curvature operators $\mathcal{R},  \mathcal{R}(g)$,
$ \mathcal{R}^* $ of the conjugate triple
$(\nabla , g, \nabla ^*)$ with torsion free connections
$\nabla , \nabla ^*$ we have  the following relations, stated  in a
standard 
local notation,  using the tensor $C$ defined in Theorem \ref{thm-7.3}. We shall omit the proof as it is 
routine \cite{EI}.
\begin{lemma}\label{lem-7.6}
We have the following identities:\begin{enumerate}
\smallbreak\item $R_{jkl}{}^i -  R(g)_{jkl}{}^i = \nabla (g)_k C_{jl}{}^i -  \nabla (g)_l C_{jk}{}^i +
C_{jl}{}^h C_{hk}{}^i - C_{jk}{}^h C_{hl}{}^i$.
\smallbreak\item $R^*{}_{jkl}{}^i -  R(g)_{jkl}{}^i = \nabla(g)_l C_{jk}{}^i -  \nabla(g)_k C_{jl}{}^i +
C_{jl}{}^h C_{hk}{}^i - C_{jk}{}^h C_{hl}{}^i$.
\smallbreak\item $
R_{jkl}{}^i - R^*_{jkl}{}^i = 2 [ \nabla(g)_k C_{jl}{}^i -  \nabla(g)_l C_{jk}{}^i]$.
\item $\mathcal{R} = \mathcal{R}^*$  if and only if  $\nabla(g) C$ 
  is 
totally symmetric.
\item $R_{jkl}{}^i +  R^{*}_{jkl}{}^i -  2 R(g)_{jkl}{}^i = 2(C_{jl}{}^h C_{hk}{}^i - C_{jk}{}^h 
C_{hl}{}^i )$.
\end{enumerate}\end{lemma}

\begin{remark}\label{rmk-7.7}\rm The generalized scalar curvature $\tau$ appears in several of
the components $\pi_j$ and $\alpha_j$. The foregoing relation
in (5) allows to calculate the deviation of $\tau$ from the Riemannian
scalar curvature $n(n-1)\kappa := \Tr_g Ric(g)$ as follows:
$$2(\tau - n(n-1)\kappa) = \|C\|^2 \, - \, n^2 \,\|T\|^2\,.$$ 
$J:=\tfrac{1}{n(n-1)}\cdot \|C\|^2$ is called the {\it Pick invariant.}  
The relation generalizes the so called {\it theorema egregium} of
relative hypersurface theory \cite{SSV}.\end{remark}

We conclude this section with:

\begin{observation}\label{obs-7.8}
 Let $(\nabla,g,\nabla{}^*)$ be a conjugate triple of torsion free connections. As $\gamma$ from \ref{rmk-7.4} and also 
$\mathcal{R}(g)$ 
are  algebraic curvature operators, we see that 
$\mathcal{R}+\mathcal{R}^*$ is a $g-$algebraic curvature operator.
\end{observation}
%%%%%%%%%%%%%%%%%%%%%%%%%%%%%%%%%%%%%%%%%%%
%%%%%%%%%%%%%%%%%%%%%%%%%%%%%%%%%%%%%%%%%%%
\section{Codazzi Structures on Manifolds}\label{sect-8}
\subsection{Conformal and projective classes}\label{sect-8.1}
Consider a pseudo-Riemannian metric $g$ and  a connection $\nabla ^*$ that is torsion free. The metric 
generates a conformal structure
 $\mathfrak{C}= \{g \}$
and the connection a projective class $\mathfrak{P^*} = \{ \nabla^* \}$ of torsion free  
connections. Any  positive function $q \in C^{\infty}(M)$  induces a simultaneous 
transformation in 
both structures, called a {\it gauge transformation}  with {\it transition function} $q$,  
by 
\begin{enumerate}
\item[(a)] a conformal change $g^\sharp = q \cdot g$;
\item[(b)] a projective change 
$$ \nabla ^{*\sharp}_v w - \nabla ^*_v w = (d \ln q)(v)w + (d \ln q)(w)v. $$
\end{enumerate}
This has the following consequences \cite{BOK-G-S,SIM-1,SIM-2,SIM-3,SIM-4}:
\begin{observation}\label{obs-8.1}
\rm\
Under the given assumptions for $\mathfrak{C}$ and $\mathfrak{P}^*$ we
have:
 \begin{enumerate} 
\item Transform a given  pair $(\nabla ^*, g)$  simultaneously 
in the conformal 
and the projective class according to the above transformations, 
respectively; 
then $(\nabla ^*, g)$ is a Codazzi pair if and only if 
$(\nabla ^{*\sharp}, g^\sharp)$ is a Codazzi pair \cite{SIM-1,SSV}.
\item
For a Ricci symmetric connection $\nabla^*$, under a  projective change with transition 
function $q$, $\nabla ^{*\sharp}$ is  
again Ricci symmetric \cite{SIM-1}. 
\item Starting with a Codazzi pair  $(\nabla ^*, g)$, 
the foregoing 
extension to $(\nabla, g, \nabla ^{*} )$ 
and transformation to 
$(\nabla^{\sharp}, g{^\sharp}, \nabla ^{*\sharp} )$ give a 
conjugate triple  with torsion free connections s.t. $(\nabla ^{*\sharp}, g^{\sharp})$ and 
$(\nabla^{\sharp}, g^{\sharp})$ are Codazzi pairs.
If additionally $\nabla^*$ is Ricci symmetric then $\nabla, \nabla^{\sharp}
\text{and} \nabla^{*\sharp}$ are Ricci symmetric.
\item In \cite{BOK-G-S} we proved that the 
foregoing transformation of conjugate triples with Codazzi pairs 
$(\nabla ^*, g)$ and $(\nabla, g)$ 
is equivalent
to a gauge transformation in an appropriate Weyl geometry. 
For this reason, for the simultaneous conformal and projective
transformations with the same transition function
within the classes $\mathfrak{C}$ and $\mathfrak{P}^*$, 
we adopt the terminology {\it gauge transformations}; invariants under
gauge transformations are called {\it gauge invariants} \cite{SIM-4}.
\item As above, consider the  conformal structure
 $ \mathfrak{C} = \{ g \} $ and the  
{\it Ricci symmetric projective structure}  $\mathfrak{P^*} = \{ \nabla^* \} $, that means the generating connection $ \nabla^*$ is torsion free and 
Ricci symmetric. If
$(\nabla ^*, g)$ is a Codazzi pair it generates a conjugate
triple $(\nabla, g, \nabla ^*)$ with Ricci symmetric connections, and
Ricci symmetric conjugate
triples go to Ricci symmetric conjugate
triples under gauge transformations. 
As in \cite{BOK-G-S}
we call  a structure, given by  a conformal
and such a projective class $\mathfrak{P^*}$, both  related by Codazzi equations, 
a {\it Codazzi structure}. The Codazzi pairing induces 
a bijective mapping $\mathfrak{P^*} \longleftrightarrow \mathfrak{C}$.
\end{enumerate}\end{observation}
%%%%%%%%%%%%%%%%%%%%%%%%%%%%%%%%%%%%%%%%%%%%%%%%%%%%55
\subsection{The  gauge invariant difference tensor}\label{sect-8.2}
In \cite{SIM-4} we listed 
 gauge invariants of conjugate triples. Above we defined
the difference tensor $ C$ and its Tchebychev form $T^{\flat};$
then 
the trace free part $\tilde{C} $  of $C$
  is a gauge invariant:
$$ \tilde{C}_{jk}{}^i := C_{jk}{}^i -    \tfrac {n}{n+2}(T_j
\delta_k^i   +  T_k \delta_j^i + g_{jk}T^i). 
$$
We use $ \tilde{C}$ in section 11.2 below.
%%%%%%%%%%%%%%%%%%%%%%%%%%%%%%%%%%%%%%%%%%%%%%%%%
\subsection{Blaschke structures}\label{sect-8.3}
Let a Codazzi structure be given by
a conformal class $\mathfrak{C}= \{g \}$
and a projective, Ricci symmetric  class $\mathfrak{P^*}= \{ \nabla^* \}$, related by 
Codazzi equations. Then there exists  a unique Codazzi pair
$(\nabla^*, g) \in \mathfrak{P} \times \mathfrak{C}$,  satisfying $\omega(g) = \omega^*$  ({\it apolarity}), where the equality holds modulo a 
non-zero constant factor. 
We call the associated conjugate triple a 
{\it Blaschke structure} or {\it equiaffine structure}
 on $M$ and use a notational mark ``e''  for {\it equiaffine}. 
The   existence
of a  Blaschke structure
follows in analogy to Proposition
5.3.1.1 in \cite{SSV}. As  $\omega(g) = c \cdot \omega^*$ with $0 \ne
c \in \mathbb{R}$ 
is equivalent to $T^{\flat} = 0$ we conclude that 
$\tilde{C} = C(e),$ that is: {\it $ C(e)$ is a trace free, gauge invariant
tensor field.}
%%%%%%%%%%%%%%%%%%%%%%%%%%%%%
\subsection{Codazzi structures and curvature operators}\label{sect-8.4}
For a conjugate triple, the connections
$\nabla$  and $\nabla ^*$
induce  curvature operators $\mathcal{R}:= \mathcal{R}(\nabla)$ and $ \mathcal{R}^*:= 
\mathcal{R}(\nabla^*)$, respectively.
$\mathcal{R}$ is equiaffine
if and only if  $\mathcal{R} ^*$ is equiaffine, see Theorem 4.10. 
Assume that both operators satisfy this condition, then both curvature
operators  $\mathcal{R} $ and  $\mathcal{R} ^*$ are conjugate, and
then a 
gauge transformation transforms a conjugate triple to a conjugate triple and thus induces a 
transformation of conjugate
curvature operators, and if $\mathcal{R}$ is equiaffine then it easily
follows from the foregoing that equiaffine  conjugate
curvature tensors 
${R}$ and $ {R}^*$
give, via gauge transformations, equiaffine  conjugate
curvature tensors ${R}^{\sharp}$ and $ {R}^{\sharp *}$.

\subsection{Equiaffine Einstein spaces}\label{sect-8.5}
 Let $(M,g)$ be pseudo-Riemannian
and $\nabla$ be a torsion free and Ricci symmetric connection. 
 $(M,g, \nabla)$ is called {\it equiaffine Einstein} if $Ric:=Ric(\nabla)$
satisfies $Ric = \lambda \cdot g$ for $\lambda \in C^{\infty}(M)$;
 this relation then holds for any other metric $g^{\sharp}$
in the conformal class of $g$.
One has the following result which
gives an equivalent condition in terms of the $W$-decomposition of Section \ref{sect-4}:
\begin{lemma}\label{lem-8.2} 
$(M,g, \nabla)$  is equiaffine Einstein if and only if 
$\pi_2({R}) = 0 = \pi_3({R})$.
\end{lemma}

\subsection{W-Decomposition,\,  A-Decomposition,
and Codazzi structures}\label{sect-8.6}
Let
$\nabla$ be a torsion free connection on a pseudo-Riemannian manifold $(M,g)$. We extend $(\nabla,g)$ to a
conjugate  triple
$(\nabla,g,\nabla^*)$;
from above we know that $\nabla^*$ is torsion free if and only
if $(\nabla,g)$  is a Codazzi pair. 
{\it From now on
we assume that both connections $\nabla,\nabla^*$ in the conjugate triple 
$(\nabla,g,\nabla^*)$ are torsion free.} 

Again, the metric $g$ generates a conformal structure $\mathfrak{C}:= \{g \},$ and the 
connection $\nabla^*$ a projective class  $\mathfrak{P^*}:= \{ \nabla^* \}$
of torsion free 
connections.

At each point $p \in M$ the tangent space is a vector space with a {scalar product}
$g_p$. The metric $g$  and both connections induce curvature operators   $\mathcal{R},$ 
 $ \mathcal{R}^*$. As $\mathfrak{P^*}$ is a projective class the projective (1,3) 
curvature operator $\mathcal{P}(\mathcal{R}^*)$ is an invariant of the class, denoted by
$$\mathcal{P}(\mathfrak{P}^*):= \mathcal{P}(\mathcal{R}^*)\quad\text{for}\quad
\nabla^* \in \mathfrak{P}^*\  \quad\text{and}\quad \mathcal{R}^* =  \mathcal{R}(\nabla^*).$$

Now we study relations between the foregoing geometric structures
on $M$ and the pointwise algebraic  W-decomposition of 
generalized curvature tensors.
Recall that  a connection $\nabla$ with curvature operator $\mathcal{R}$
is called {\it Ricci symmetric}  if its Ricci tensor $Ric$ is symmetric at any point  $p 
\in M$; as already stated, this is equivalent to the fact that $\nabla$ locally admits a 
parallel volume form.
\begin{theorem}\label{thm-8.3}
 A conjugate triple $(\nabla,g,\nabla^*)$, its
induced curvature operators $\mathcal{R}$ and $\mathcal{R}^* $,  and their
decompositions  satisfy the following equivalences:
\begin{enumerate}
\item  $\nabla$ is Ricci symmetric. 
\item $\nabla$ locally admits a parallel volume form 
$\omega$, thus $\nabla \omega = 0$. 
\item $R$ is an equiaffine curvature tensor.
\item  $\nabla^*$ is Ricci symmetric. 
\item $\nabla^*$ locally admits a parallel volume form 
$\omega^*$, thus $ \nabla^*  \omega^* =0$.
\item $R^*$ is an equiaffine curvature tensor. 
\item $\pi_3({R})=0= \pi_8({R})$ at any point.
\item $\pi_3({R}^*)=0= \pi_8({R}^*)$ at any point.
\item $\alpha_4({R})=0=\alpha_5({R})=\alpha_8({R})$ at any point.
\item $\alpha_4({R}^*)=0=\alpha_5({R}^*)=\alpha_8({R}^*)$ at any point.
\end{enumerate}
\end{theorem}
\begin{proof}
It follows from Theorem  \ref{thm-7.3}     and  Remark   \ref{rmk-7.4} that both curvature tensors $R$ and $R^*$ satisfy
$R, R^* \in \mathfrak{r}(T_pM)$ on $M$. Theorem \ref{thm-4.10}
implies that $\nabla$ is Ricci symmetric if and only if 
$\nabla^*$ is Ricci symmetric \cite{SSV}.
For the rest of the proof see the  results given above.\end{proof}

\begin{remark}\label{rmk-8.4}
\rm The relation $\pi_3({R})=0$ implies  $\pi_4({R})=0$; analogously 
$\pi_3({R}^*)=0$ implies  $\pi_4({R}^*)=0.$
\end{remark}
%%%%%%%%%%%%%%%%%%%%%%%%%%%%%%%%%%%%%%%%%%%%%%%%%%%%%%%%%%%%%%%
\section{Projective and Conformal Changes}\label{sect-9}
%%%%%%%%%%%%%%%%%%%%%%%%%%%%%%

On a given pseudo-Riemannian manifold we consider a torsion free connection $\nabla^*$ and 
the projective class
$\mathfrak{P}^* $ of torsion free connections generated
by $\nabla^*$. As before  we write $\mathcal{R}^* := \mathcal{R}(\nabla^*)$ etc.

\subsection{\bf W-Decomposition and projective classes}\label{sect-9.1}
From  \cite{BO} we have:
\begin{theorem}\label{thm-9.1}
Consider a projective class $\mathfrak{P}^*$ of torsion free 
connections; for any 
$\nabla^* \in \mathfrak{P}^* $ the five 
components $\pi_j({R^*}),$ for $j=4,...,8,$
in the $W$-decomposition of the  (0,4) 
curvature tensor 
$$ {p}(\mathfrak{P}^*) =  {p}({R}^*) = \sum_{j=4}^8 \pi_j({R}^*)$$
give rise to projective invariants, namely their $g$-associated (1,3) curvature operators  in 
the associated decomposition
of the projective (1,3) curvature  operator in $\mathfrak{R}(V).$
Analogously,
the (1,3) curvature   operator that is $g$-associated to 
$${R^*}  -  \left(\pi_1({R}^*) \oplus \pi_2({R}^*) \oplus 
\pi_3({R}^*)\right),$$
is a projective invariant.
\end{theorem}
\begin{proof} For all $\nabla^* \in  \mathfrak{P}^*$ the projective curvature
operators  coincide, and according to Section \ref{sect-2} one
can characterize the subspaces in the  $W$-decomposition,
and from this  the components $\pi_j({R}^*)$ for $j=4,...,8$  are uniquely
determined,  and  we have
$$ \pi_j({R}^*) = \pi_j({p}({R}^*))\,.$$
Thus the components in the  associated  decomposition 
of the projective (1,3) curvature tensor
are the same for any 
$\nabla^*$, that means that their $g$-associated (1,3) curvature operators 
are projectively invariant.
 \end{proof}
{\bf Remark.} \, In case that we study equiaffine curvature tensors, we
have 
$$ {p}(\mathfrak{P}^*) =  p({R}^*) = \sum_{j=5}^7 \pi_j(R^*),$$
and each component and also 
$${R^*}  -  \left(\pi_1({R}^*) \oplus \pi_2({R}^*)\right)$$
lead to   projective invariants taking the 
the associated decomposition
of the projective (1,3) curvature operator in $\mathfrak{F}(V)$.

The following Theorem is a corollary of
the foregoing decomposition in the space  $\mathfrak{F}(V),$   but because of its importance
it is stated as a separate result; it generalizes results of \cite{STEG}. 
Here  we get a projectively  invariant symmetric bilinear form
of a Codazzi structure in a very general situation. See the material above in Section \ref{sect-4}.
We recall  definition \ref{defn-4.13}  and that  of $B^*$ from Equation (\ref{eqn-4.b}).

\begin{theorem}\label{thm-9.2} The symmetric bilinear form 
$B^*$
is invariant under a projective change in    the Ricci symmetric class $\mathfrak{P}^*$.
\end{theorem}

\begin{proof} The (1,3) curvature operator  in $\mathfrak{F}(V)$ 
that is $g$-associated to $\pi_5({R^*}) \subset \mathfrak{f}(V,g)$
is a projective invariant. Take the trace $\rho_{13}$.\end{proof}

The following results are useful as well; as before,  in this  section
 we assume 
that the connections $\nabla, \nabla^*$ are torsion free:
\begin{lemma}\label{lem-9.3}
Let $(\nabla,g,\nabla^*)$ be a conjugate triple
with Ricci symmetric $\nabla^*$. Then
$$p({R}^*) = \sum_{j=5}^7 \pi_j({R}^*);$$
if additionally $\mathcal{P}(\mathcal{R}) = \mathcal{P}(\mathcal{R}^*)$ then  $Ric = Ric^*$.
\end{lemma}

The foregoing statements  give  the following result:
\begin{corollary}\label{cor-9.4}
 Let $\nabla^*$ be Ricci symmetric and  $(\nabla,g,\nabla^*)$
be a conjugate triple as above. The following equations are equivalent:
\begin{enumerate}\item $\mathcal{P}(\mathcal{R}) = \mathcal{P}(\mathcal{R}^*)$.
\item $\mathcal{R} = \mathcal{R}^*$.
\end{enumerate}
\end{corollary}

We now discuss the $W$-Decomposition   and projective changes.
According to Weyl two connections $\nabla^*$ and $\nabla^{*\sharp}$
are projectively equivalent if and only if there exists a one-form
$\theta$ s.t.
$$\nabla^*_vw - \nabla^{*\sharp}_{\,v}\,w  =  \theta(v)w +
\theta(w)v\,.$$
In section \ref{sect-8.1} we considered the special case of $\theta =
d\ln\,q,$  where $q \in C^{\infty}$ is a  transition function.
Observation 8.1 gives  the following:

\begin{proposition}\label{prop-9.5}
Let $(\nabla, g, \nabla^*)$ be a conjugate triple
and  $0 < q \in C^{\infty}$ a transition function  so that
$(\nabla, g, \nabla^*) \mapsto (\nabla^{\sharp}, g^{\sharp}, \nabla^{*\sharp})$
is a Codazzi transformation. The following equations are equivalent:
\begin{enumerate}
\item $\pi_3({R}) = 0$.
\item $\pi_3({R^*}) = 0$.
\item $\pi_3({R}^{\sharp}) = 0$.
\item $\pi_3({R}^{\sharp*}) = 0$.
\end{enumerate}\end{proposition}
We relate the $W$-Decomposition and projective flatness. Recall that a projective class $\mathfrak{P}^*:=
\{\nabla^*\}$  that is generated from  $\nabla^*$ by gauge transformations
with  transition functions is said to be  a {\it Ricci symmetric projective class}
if $\nabla^*$ is torsion free and Ricci symmetric. From the foregoing
then  {\it any}  connection in   $\mathfrak{P}^*$  is  torsion free and Ricci symmetric.
It is well known that  projective flatness can be 
characterized   as follows \cite{EIS}:
\begin{lemma}\label{lem-9.6}
 Let $n \ge 2$ and $\nabla^*$ be
torsion free and Ricci symmetric. Then
$\nabla^*$ is {\it projectively flat}  if and only if the projective curvature 
tensor $ \mathcal{P}(\mathcal{R^*})$ vanishes and the covariant derivative
$\nabla^*\,Ric^*$ is totally symmetric
\end{lemma}

It is also well known that the two conditions in Lemma \ref{lem-9.6} are dependent; see e.g. \cite{EIS,NOM-SA,SCHI,SSV}:
\begin{lemma}\label{lem-9.7}
\begin{enumerate}
\item In dimension $n=2$ the projective curvature tensor
vanishes identically and projective flatness is equivalent to
the total symmetry of $\nabla^*\,Ric^*$.
\item In dimension $n \ge 3$ projective flatness is equivalent to
the vanishing of the projective curvature tensor $ \mathcal{P}(\mathcal{R^*})$; the total symmetry of 
$\nabla^*\,Ric^*$ is a consequence.
\end{enumerate}\end{lemma}
\begin{lemma}\label{lem-9.8} Let  $n \ge 3$ and $\nabla^*$ be Ricci symmetric. Then the
following assertions are equivalent:
$(1)\ \nabla^*$ is projectively flat;  $(2)\ \mathcal{P}(\mathcal{R^*}) = 0$;
$(3)\ {R}^* \in W_1 \oplus W_2$.
\end{lemma}
The following Corollary
is now immediate; it is a technical use we shall need subsequently.

\begin{corollary}\label{cor-9.9}
Let $n \ge 3$ and $(\nabla, g, \nabla^*)$ be a conjugate triple; 
assume that $\nabla^*$ is  Ricci symmetric and projectively flat. Then:
$$ B^*:= (n-1) Ric + Ric^* - \tau g=0\,.$$
\end{corollary}

We now have:
\begin{theorem}\label{thm-9.10}
{\bf[Equivalence Theorem]}
 Let $n \ge 3$ and $(\nabla, g, \nabla^*)$ be a conjugate 
triple; 
assume that $\nabla^*$ is  Ricci symmetric and projectively flat. Then
the following assertions are equivalent:
\begin{enumerate}
\item $\nabla$ is projectively flat.
\item $Ric = Ric^*$.
\item $\pi_2({R}) =0$.
\item $\pi_2({R}^*) =0$.
\item $n\cdot Ric  \,= \tau \cdot g$; that means: $(\nabla,g)$ is equiaffine Einstein.  
\item $n\cdot Ric^* = \tau \cdot g$; that means: $(\nabla^*,g)$ is equiaffine Einstein.  
\item $\mathcal{R}  = \mathcal{R}^*$.
\item $-\, n(n-1)\cdot {R} = \tau \cdot (g \wedgeo  g)$.
\item The covariant derivative $\nabla(g)C^{\flat}$ is a totally symmetric 
(0,4)-tensor field.
\end{enumerate}\end{theorem}

We can draw the following consequence:

\begin{corollary}\label{cor-9.11} 
Let $n \ge 3$ and $(\nabla, g, \nabla^*)$ be a conjugate triple; 
assume that $\nabla^*$ is  torsion free, Ricci symmetric, projectively flat and equiaffine Einstein. Then $\tau =
const.$
\end{corollary}
\begin{proof} $\nabla^*$ is  torsion free, thus $(\nabla, g)$ and also 
$(\nabla^*, g)$ are Codazzi pairs. The projective flatness implies
that also $(\nabla^*, Ric^*)$ is a Codazzi pair \cite{EIS}.
$\nabla^*$-covariant differentiation of the equation $Ric^* = \tfrac{1}{n}\cdot \tau \, 
g$,  and 
the Codazzi properties  together with a 
contraction finally give the assertion.
\end{proof}

The following result  is another  simple consequence of Theorem \ref{thm-9.10};
because of its geometric importance we state it as a Theorem; namely, it is remarkable, 
that, under the given assumptions, the projective flatness of $\nabla$
is equivalent to the vanishing of a symmetric bilinear form, and there is no need
to calculate its  (1,3) projective curvature operator.

\begin{theorem}\label{thm-9.12} Let $n \ge 3$ and $(\nabla, g, \nabla^*)$ be a conjugate triple; 
assume that $\nabla^*$ is  torsion free, Ricci symmetric and  projectively flat.
Then the following conditions are equivalent:
\begin{enumerate}
\item  $\nabla$ is projectively flat.
\item $B:= (n-1)Ric^* + Ric - \tau\cdot g = 0.$
\end{enumerate}\end{theorem}

\begin{proof}The projective flatness of $\nabla^*$ implies $B^* = 0$. Then $B = 0$ is equivalent
to the identity $Ric^* = Ric$. An easy computation now completes the proof. \end{proof}

%%%%%%%%%%%%%%%%%%%%%%%%%%
\subsection{Projective flatness and PDEs}\label{sect-9.2}
 Theorem \ref{thm-9.10}\,(9) relates properties of equi\-affine curvature tensors with
that of PDEs. This admits important applications. We give the following example  that 
generalizes the local classification of locally strongly convex equiaffine
spheres with constant sectional curvature of the Blaschke metric. It is
remarkable 
 that, in the context of conjugate connections, the essential assumptions can be 
expressed in terms of the $W$-decomposition of the
three curvature operators $\mathcal{R}$, \, $\mathcal{R}(g)$, and 
$\mathcal{R}^*.$ The restriction to dimension $n \ge 3$ is due to the fact that in the 
following we characterize projective flatness by
the  vanishing of the projective curvature tensor.
\begin{theorem}\label{thm-9.13} Let $n \ge 3$ and $(\nabla, g, \nabla^*)$ be a Blaschke structure
with a Riemannian metric $g$.  Assume that ${R}(g) = \pi_1({R}(g))$, that ${R}^* = \pi_1({R}^*)
\oplus \pi_2({R}^*)$, and that $\pi_2({R}) = 0$.
Then: 
\begin{enumerate}
\item If $C = 0$ then trivially all three connections coincide.
\item If $\|C\|  \ne 0$ then $(M,g)$ is flat and $\tau$ is a  negative constant. 
\end{enumerate}
\end{theorem}
\begin{proof}
 From the assumptions, from the results of Section \ref{sect-7}  
and from Theorem \ref{thm-9.10}
the following conditions are satisfied:
\begin{enumerate}
\item $C^{\flat}$ is totally symmetric.
\item $\Tr(C) = 0$.
\item $g(C(u,v),C(w,z)) \, - \, g(C(w,v),C(u,z))$\smallbreak$=
  ({\tau} - \kappa)(g(u,v)\,g(w,z) - g(w,v)\,g(u,z))$.
\item $\tau = const.$ and $\kappa = const$ (see Remark 7.7).
\item $(\nabla(g)_{u} C)(v,w) = (\nabla(g)_{v} C)(u,w)$.
\end{enumerate}
The proof follows now the lines of the proof of the Main Theorem in \cite{VLS}.
\end{proof}

Concerning the local classification of equiaffine spheres with indefinite metric and constant 
Blaschke sectional curvature, there is the famous solution of the so called Magid-Ryan 
conjecture by Vrancken \cite{VRA-LN}, \cite{VRA}.
This result and its proof can be generalized to conjugate connections
as follows:

\begin{theorem}\label{thm-9.14} Let $n \ge 3$ and $(\nabla, g, \nabla^*)$ be a 
Blaschke structure with indefinite metric $g$. Assume that
\begin{enumerate}
\item $\pi_2({R}) = 0$.
\item ${R}(g) = \pi_1({R}(g))$.
\item ${R}^* = \pi_1({R}^*) \oplus \pi_2({R}^*)$.
\end{enumerate}
Moreover, assume that $\tau - \kappa \ne 0.$ Then $(M,g)$ is flat.
\end{theorem}

\begin{proof}As in the proof of Theorem \ref{thm-9.13}, the conditions
(1) - (5) are satisfied. Apply now Theorem 6 in \cite{VRA-LN}.\end{proof}

\begin{theorem}\label{thm-9.15} Let $n \ge 3$ and $(\nabla, g, \nabla^*)$ be a 
Blaschke structure with indefinite metric $g$. Assume that
$\pi_5({R}) = 0$, assume that ${R}(g) = \pi_1({R}(g))$, and assume that
${R}^* = \pi_1({R}^*) \oplus \pi_2({R}^*)$.
 Then $(M,g)$ is flat and $\|C\|^2 = 0$.
\end{theorem}

\begin{proof}
Apply Theorem 11 in \cite{VRA-LN}.
Analogously one can generalize Theorem 12 in Vrancken's paper \cite{VRA-LN}.
\end{proof}

\subsection{W-Decomposition and gauge transformations of conjugate triples}\label{sect-9.3} 
We studied ``pointwise'' conformal changes in Section \ref{sect-3.2}, while 
projective changes on a manifold  where investigated above.
It is well known that the Ricci symmetry is invariant under a projective change with a 
transition function \cite{PSS}. 
Considering a conjugate triple $(\nabla,g, \nabla^*)$
 and their (1,3) curvature operators  $\mathcal{R}$ and $ \mathcal{R}^*,$ and
also the curvature operator $\mathcal{R}(g)$ of the metric $g$, Section \ref{sect-8.1} and the foregoing 
results
give the following Theorem.

\begin{theorem}\label{thm-9.16}
Let $\nabla$ and $ \nabla^*$ be torsion free connections.
Let $$(\nabla,g, \nabla^*) \mapsto (\nabla^{\sharp},g^{\sharp}, 
\nabla^{*\sharp})$$
be a gauge transformation with a transition function as in Section \ref{sect-8.1}. Then:
\begin{enumerate}
\item If $\nabla$  or $\nabla^*,$ respectively, is Ricci symmetric,
then all connections appearing in the conjugate triple and under gauge
transformations of this triple are Ricci symmetric.
\item We have that
\medbreak\centerline{
$\displaystyle\sum_{j=4}^7 \pi_j({R}^*)  = p({R}^*)  =  p({R}^{*\sharp}) 
= 
 \sum_{j=4}^7 \pi_j({R}^{*\sharp})$.}\medbreak\noindent
In particular: each component of the decompositions  of (1,3) curvature
operators  that is $g$-associated to $\pi_j({R}^*)$,  for
$j = 4,...,7$,   is gauge 
invariant itself;
if $\nabla$ or  $\nabla^*$, respectively, is Ricci symmetric  then 
$\pi_4({R}) = 0 = \pi_4({R}^*)$.
\item The conformal class satisfies:
$$ \alpha_6({R}(h))=  \pi_6({R}(h)) = 
\pi_6({R}(h^{\sharp}))=  \alpha_6({R}(h^{\sharp}))\,.$$
\end{enumerate}
\end{theorem}

\begin{corollary}\label{cor-9.17}
Let $\nabla$ and $\nabla^*$ be torsion free and Ricci symmetric
connections. 
Let $$(\nabla,g, \nabla^*) \mapsto
(\nabla^{\sharp},g^{\sharp}, 
\nabla^{*\sharp})$$
be a gauge transformation with a  transition function. If for one
of the curvature tensors, say ${R}$, we have  
$\pi_j({R}) = 0$ for $j = 3 $ \, ($j = 4 $, resp.)  then for any
curvature tensor under conjugation or gauge transformation, the
components $\pi_j$ vanish for $j =3$ \, ($j = 4 $, resp.).
\end{corollary}

We also have:
\begin{theorem}\label{thm-9.18} Let   $(\nabla,g, \nabla^*)$ be a
  conjugate triple with $\nabla^*$ torsion free and Ricci symmetric. Then the 
bilinear form $B^*$ in Corollary \ref{cor-9.9} is  gauge invariant.
\end{theorem}

\begin{proof}For fixed metric $g$ the (1,3) curvature operator, that is 
$g$-associated to 
the component $\pi_5({R}^*)$,  is projectively 
invariant. Following (2) in the foregoing Theorem, a conformal 
change $g^{\sharp} = \lambda \cdot g$ gives  $\pi_5^{\sharp}({R}^*) =  
\pi_5({R}^*)$, thus the associated (1,3) curvature operators 
in $\mathfrak{R}(V)$ are independent of the conformal change. Taking traces
$\rho_{13}$ on both sides gives $B^{*\sharp} = B^*$. 
\end{proof}
%%%%%%%%%%%%%%%%%%%%%
\subsection{The decomposition of equiaffine curvature tensor
  fields}\label{sect-9.4}
According to its geometric importance, we summarize  the following
observations in the special case  of equiaffine curvature tensor
fields.  Throughout we assume  that $(\nabla,g, \nabla^*)$ be a conjugate triple with torsion free 
connections $\nabla$ and $\nabla^*$.
\begin{observation}\label{obs-9.19}
 If $\nabla^*$ is Ricci symmetric then 
$\mathcal{R}^*$ is an equiaffine curvature operator, and  then
$\mathcal{R}$ is an equiaffine curvature operator.
\end{observation}
For the $W$-decomposition one has that:
\begin{observation}\label{obs-9.20}
\ \begin{enumerate}
\item ${R}$ and  ${R}^*$, respectively,  satisfy the orthogonal
  decompositions
\begin{eqnarray*}
&&{R} = \pi_1({R})   \oplus   \left(\pi_2({R}) \oplus 
\pi_5({R})\right)   \oplus   \left(\pi_6({R}) \oplus \pi_7({R})\right), \\
&&{R}^* = \pi_1({R}^*)   \oplus  \left(\pi_2({R^*}) \oplus 
\pi_5({R^*})\right)  \oplus  
\left(\pi_6({R}^*) \oplus \pi_7({R}^*)\right).
\end{eqnarray*}
\item The sums   $(\pi_2({R})  \oplus  \pi_5({R}))$   and   
$(\pi_2({R}^*)  \oplus \pi_5({R}^*))$   
are   Ricci traceless and also Ricci${}^*$ traceless. 
\item The sums  
$ \left(\pi_6({R}) \oplus \pi_7({R})\right)$ 
  and 
$\left(\pi_6({R}^*) \oplus \pi_7({R}^*)\right)$  
are Ricci flat and also Ricci${}^*$ flat.
\end{enumerate}
\end{observation}
For the $A$-decomposition, we have:
\begin{observation}\label{obs-9.21} 
${R}$ and  ${R}^*$, respectively,  satisfy the orthogonal decomposition
\medbreak\qquad
$\displaystyle {R} =  \alpha_1({R})   \oplus    (\alpha_2({R}) \oplus 
\alpha_3({R}) )   \oplus   (\alpha_6({R}) \oplus \alpha_7({R}) ), $
\medbreak\qquad
$\displaystyle {R}^* =  \alpha_1({R}^*)   \oplus    (\alpha_2({R}^*) \oplus 
\alpha_3({R}^*) )   \oplus    (\alpha_6({R}^*) \oplus \alpha_7({R}^*) )$
\medbreak\noindent 
where
\begin{enumerate}
\item both of the two orthogonal sums $\alpha_2({R}) \oplus \alpha_3({R})$ and 
$\alpha_2({R}^*) \oplus \alpha_3({R}^*)$ are 
Ricci traceless and Ricci${}^*$ traceless;
\item both of  the two orthogonal sums
$ \alpha_6({R}) \oplus \alpha_7({R}) $ and 
 $ \alpha_6({R}^*) \oplus \alpha_7({R}^*)$
are  Ricci flat and Ricci${}^*$ flat;
\item $\pi_1({R}^*) = \pi_1({R}) = 
\alpha_1({R}) = \alpha_1({R}^*)$  
and  $\pi_6({R}) = \pi_6({R}^*)$;
\item $ \pi_2({R}) \oplus  \pi_5({R}) = 
\alpha_2({R}) \oplus  \alpha_3({R})$
\newline\qquad$=\frac{1}{n(n-1)}[2\tau\, 
g\wedgeo g +  g\wedge_{n-1}Ric-Ric^*\wedge_{n-1}g]$;
\item $ \alpha_5({R}) = \alpha_5({R}^*)$.
\end{enumerate}\end{observation}
Concerning projective invariants, Theorem \ref{thm-9.1} 
now reads:
\begin{observation}\label{obs-9.22}
The (0,4) curvature tensor  $p({R}^*)$  satisfies
$$p({R}^*) = \sum_{5}^7 \pi_j({R}^*)\,.$$
The (1,3) tensor components $g$-associated to $\pi_j({R}^*)$ are projective 
invariants for $j = 5,6,7$.
\end{observation}
\begin{remark}\label{rmk-9.23}\rm We would like to comment on the foregoing summary of
decomposition results; for this, we recall the decomposition of algebraic curvature tensors 
in section 5.1 and that of generalized curvature tensors discussed in Section \ref{sect-6}.  Considering conjugate
connections and their equiaffine curvature tensor  fields on a manifold,  we see how the decomposition reflects
geometric properties of a triple $(\nabla,g, \nabla^*)$; moreover, we learn that the concepts
of conjugate connections and conjugate curvature tensors, the latter induced
from the first, are appropriate tools to understand how properties
of algebraic curvature tensors generalize to equiaffine curvature tensors.
\end{remark}
\newpage
\section{Relative Hypersurface Theory}\label{sect-10}
We recall basics from relative hypersurface theory \cite{LI-S-Z-1,SSV}.
%%%%%%%%%%%%%%%%%%%%
\subsection{Review of  relative hypersurface theory}\label{sect-10.1}
We  describe 
the duality of the vector  space $ \mathbb{R}^{n+1}$ and its
dual $\mathbb{R}^{(n+1)*}$  in terms of a
non-degenerate scalar product
$$
\langle \quad  , \quad \rangle: \mathbb{R}^{(n+1)*} \times
\mathbb{R}^{n+1} \rightarrow \mathbb{R}.
$$
Associated to each of the vector spaces there is a one-dimensional
vector space of determinant forms, fixing volumes modulo scaling.
By  $det$ and $det^*$ we denote an arbitrary  pair of dual
determinant forms on $\mathbb{R}^{n+1}$ and $\mathbb{R}^{(n+1)*}$.
By the same symbol $\overline{\nabla}$ we denote the
canonical flat connections on $\mathbb{R}^{n+1}$ and
$\mathbb{R}^{(n+1)*}$, respectively. 
For a hypersurface immersion $x:M\rightarrow
\mathbb{R}^{n+1}$ we define a {\it normalization}: it is a pair $(Y,z)$ with
$\langle Y,z\rangle=1$ where $z:M \rightarrow   \mathbb{R}^{n+1}$ is an
arbitrary {\it transversal field}, and $Y:M \rightarrow
\mathbb{R}^{(n+1)*},$ satisfying $\langle Y,dz(v)\rangle = 0 $ for
all tangent vectors $v$ on $M$, is a {\it conormal field} of $x.$
While a transversal field $z$ extends a tangential basis to the
ambient space, a conormal fixes the tangent plane. A {\it
normalized hypersurface} is a triple $(x,Y,z).$ 
%%%%%%%%%%%%%%%%%%
\subsubsection{Structure equations}\label{sect-10.1.1}
The geometry of
$(x,Y,z)$ can be described in terms of geometric
invariants defined via the {\it structure equations}  of {\it
Gau{\ss}}  and {\it Weingarten},
 respectively:
\begin{eqnarray}
\overline{\nabla} _v dx(w) & = & dx(\nabla_v w)   + \,h(v,w)z,
\nonumber \\
dz(v) & = &  dx(-S(v))  + \sigma (v) z  \nonumber.
\end{eqnarray}
As before $u,v,w,...$ denote tangent vectors and fields,
respectively. The {\it induced connection} $\nabla$ is torsion
free, $h$ is bilinear and symmetric, $S$ is the {\it shape} or
{\it Weingarten operator} and $\sigma$ 
is a 1-form, the {\it connection form}; the sign in front of $S$ in the Weingarten
equation is a convention corresponding to an appropriate choice of
the orientation of $z$. All coefficients in the structure
equations depend on the normalization, they are invariant under
the affine group of transformations in $  \mathbb{R}^{n+1}.$
%%%%%%%%%%%%%%%%%%%%%%%%%
\subsubsection{Non-degenerate hypersurfaces}\label{sect-10.1.2}
A hypersurface $x$ is {\it non-degenerate} if the bilinear form
$h$ in the Gau{\ss} structure equation is non-degenerate; it is well
known that this property is independent of the choice of the
normalization as all such symmetric bilinear forms are conformally
related, defining a {\it conformal class} $\mathfrak C$. Thus, on
a non-degenerate  hypersurface, any transversal field induces a
pseudo-Riemannian {\it metric}  $h \in \mathfrak C$  with Levi-Civita
connection  $\nabla (h)$  and Riemannian volume form $\omega (h)$; similarly we
denote its curvature tensor by $R(h)$, its Ricci tensor by
$Ric(h)$, its normalized scalar curvature by $\kappa (h)$,  etc.

The non-degeneracy of $x$ is equivalent to the fact that any
conormal field $Y$ itself is an immersion $Y: M \rightarrow
\mathbb{R}^{(n+1)*}$ with transversal position vector $Y$. The
associated {\it Gau{\ss} structure equation}  reads
$$\overline{\nabla} _v dY(w) = dY(\nabla^*_v w) + \tfrac{1}{n-1}\,
Ric^* (v,w)(-Y)$$ where the {\it conormal connection} $\nabla^*$
is torsion free and {\it Ricci symmetric}. 
It is well known that all conormal
connections are projectively related; we denote the {\it projective class of
all conormal connections} by $\mathfrak P^*.$
%%%%%%%%%%%%%%%%%%%%%%%%%%%%%%%%%%%%%%%%%%%%%%%%%%%%%%%%%%%%%%%
%%%%%%%%%%%%%%%%%%%%%%
%%%%%%%%%%%%%%%%%%%%%%%%%
\subsubsection{Relative normalizations and curvature
  operators}\label{sect-10.1.3}
Within the class of all normalizations of a non-degenerate
hypersurface there is a distinguished large subclass, namely the
class of all {\it relative normalizations}. This class can be
characterized by the property that $\sigma = 0$ in 
the Weingarten
structure equation. This is equivalent to the fact that
the triple $(\nabla,h,\nabla^*)$ is conjugate.
In the following we restrict to this class;
this can be geometrically justified \cite{SIM-4}.
We denote a relative normalization by $(Y,y)$
and call such a triple $(x,Y,y),$ where $x$ is non-degenerate, a {\it
relative hypersurface}.
\newpage
We  recall the notation
$$\mathcal{R}:= \mathcal{R}(\nabla),\quad\mathcal{R}^*:= \mathcal{R}(\nabla^*),\quad
\mathcal{R}(h)$$ for the  
curvature operators.
Note  that, if we have an arbitrary normalization that is not relative,
then $\nabla $ is not  Ricci symmetric and thus $\mathcal{R}$
is not equiaffine, while  $\nabla^*  $ is always Ricci symmetric and
thus  $\mathcal{R}^*$ equiaffine; \cite{SIM-4}.  
Using Theorem \ref{thm-4.10}, 
  we are able 
to characterize the important  class of relative normalizations in terms
of their $W$-decomposition as follows:
\begin{lemma}\label{lem-10.1}
A normalization of a non-degenerate hypersurface
is relative if and only if $\pi_3({R}) = 0$.
\end{lemma}
%%%%%%%%%%%%%%%%%%%%%%%%%%%%%%%%%%%%%%%%%%%%%%%%%%%%%%%
\subsubsection{The cubic form and the Tchebychev form}\label{sect-10.1.4}
We recall section \ref{sect-7}.
 $$ C(v,w) = \nabla(h)_vw -  \nabla^*_vw$$
is a symmetric (1,2)-tensor  field (both connections are torsion free), and its trace 
$$nT^{\flat}(v)= \Tr(w \mapsto C(v,w))$$
is a closed 1-form as both connections locally admit parallel volume forms.
$T^{\flat}$ is called the {\it Tchebychev form}, the associated
vector field $T$, implicitly defined by $h(T,v):=T^{\flat}(v),$
 is called the {\it Tchebychev field}.
Associated to $C$ is the totally symmetric
{\it cubic form} $C^{\flat}$, defined by $C^{\flat}(u,v,w) := h(C(u,v),w)$.
%%%%%%%%%%%%%%%%%%%%%%
\subsubsection{Relative Gau\ss\  maps}\label{sect-10.1.5}
In case of a relative normalization we
know that the shape operator $S$ is $h$-selfadjoint and satisfies
$$
Ric^*(v,w) = (n-1) h(Sv,w) =:(n-1)S^{\flat}(v,w).
$$
The symmetric bilinear form $S^{\flat}$ is called the {\it
Weingarten form}. In the case that for a relative normalization
$rank\, S = n$, both,  the relative normal and conormal fields
$$y: M  \rightarrow \mathbb{R}^{(n+1)}\quad\text{and}\quad
Y:  M  \rightarrow\mathbb{R}^{(n+1)*},$$  are immersions, called {\it
Gau\ss\  maps}. For a relative hypersurface $(x,Y,y)$ the induced triple
$(\nabla,h,\nabla^*)$ is a conjugate triple with torsion free and
Ricci symmetric connections $\nabla, \, \nabla^*$.
Thus we can apply the results stated in sections \ref{sect-7} 
and \ref{sect-8}. Recall that $\Tr(S)=nH$ where $H$
denotes the {\it relative mean curvature}.
%%%%%%%%%%%%%%%%%%%%%
\subsubsection{Integrability conditions}\label{sect-10.1.6}
The integrability conditions for a relative
hypersurface 
can be stated as follows (compare  Lemma 9.6 - 9.7): 
\begin{enumerate}
\item The conormal connection $\nabla^*$ is projectively flat,
that means the (1,3) {\it projective curvature operator}
$P^*:= \mathcal{P}(\mathcal{R}^*)$  
\begin{eqnarray*}
P^*(u,v)w
 :=\mathcal{R}^*(u,v)w -\textstyle\frac1{n-1} [Ric^*(v,w)u - Ric^*(u,w)v]
\end{eqnarray*} 
 vanishes identically,  and 
 $(\nabla^*,Ric^*)$ is a   Codazzi pair, see Definition
\ref{defn-7.2}.
\item  $(\nabla^*, h)$ is a Codazzi pair.
\end{enumerate}
\subsection{Examples of relative normalizations}\label{sect-10.2}
We refer to \cite{SSV} for 
well known 
examples and further details for  relative normalizations.
%%%%%%%%%%%%%%%%%%%%%%%%%%%%%%%%%
%%%%%%%%%%%%%%%%%%%%%%%%%%%%%%%%%%%%%%%%%
\subsubsection{The equiaffine (Blaschke) normalization}\label{sect-10.2.2}
There is a (modulo sign)  unique normalization within all relative
normalizations, characterized by the vanishing of its Tchebychev field:
  $T(e) = 0$ ({\it apolarity condition});
following section \ref{sect-8.3} 
here  we use  the
notational  mark $''e''$ for the induced
{\it equiaffine} geometry. The transversal field $y=y(e)$ in
this normalization historically is called {\it
affine normal field}.
Equivalent  to the equation $T(e) = 0$ is
the relation $\omega^* = \omega (h)$  (modulo a positive constant factor).
This relation characterizes a unique Codazzi pair $(\nabla^*, h)$
within the Cartesian product  $\mathfrak P^* \times \mathfrak C.$
%\newline
Nowadays the unimodular  geometry
is often  called {\it Blaschke geometry}.
The geometry induced from the Blaschke normalization is invariant under
the unimodular transformation group (including
parallel translations).
%%%%%%%%%%%%%%%%%%%%%%%%%%%%%%%%%%
\subsubsection{The centroaffine normalization}\label{sect-10.2.3}
For a non-degenerate hypersurface  it is well known that the set
$\{p \in M  \mid  x(p)  \quad  tangential \}$ is nowhere dense. Thus
the position vector $x$ is  transversal almost everywhere.
We call a non-degenerate hypersurface $x$ with always  transversal position
vector {\it centroaffine}, and  denote the position vector
also by $x$ \cite{NOM-SA}. For such a
 hypersurface one can choose $y(c):= \varepsilon x$ as relative normal
where $ \varepsilon = +1$ or $ \varepsilon = -1$ is chosen
appropriately; 
in analogy to the foregoing we use ``$c$'' as a mark in case of a
{\it centroaffine normalization}  $(Y(c), y(c)).$  $Y(c)$ is oriented
such that
$$ 1 = \langle Y(c),y(c)\rangle. $$
%%%%%%%%%%%%%%%%%%%%%%%%%%%%%%%%%%%%%%%%%%%%%%%%%%%%%%%%%%%%%%%%%%
\subsection{Gauge invariance}\label{sect-10.3}
{}From the foregoing it is obvious that the conformal and the
projective structure are of particular importance  in relative
hypersurface theory; both classes do not depend on a particular
choice of a normalization, thus it is of interest how their
invariants reflect the geometry of a given hypersurface \cite{SIM-4}. Following the terminology of \cite{SIM-4},
{\it gauge invariants} are invariants that do not depend on a particular
choice of a normalization.
In relative hypersurface theory, the class $\mathfrak P^*$ is 
torsion free, Ricci symmetric and projectively flat; the last  
geometric
property is equivalent to one of the integrability conditions of
the structure equations, and this equivalence gives a geometric
understanding of a version of the relative fundamental theorem that is an
extension of the  original  result of Dillen, Nomizu and Vrancken
in the theory of Blaschke hypersurfaces \cite{SSV,SIM-4}. From this the projective
class  $\mathfrak P^*$   and its geometry are well understood.

The situation is different with the conformal class $\mathfrak C.$
One knows that $\mathfrak C$ is a class of Riemannian metrics if
and only if $x$ is locally strongly convex; this implies that a
connected, closed (compact without boundary) hypersurface with
definite class $\mathfrak C$ is a {\it hyperovaloid}. But e.g., so
far there is no characterization of the class of hypersurfaces for
which  $\mathfrak C$ is locally conformally flat, even  not under
strong additional assumptions like locally strong convexity and
compactness. One knows many local examples of hypersurfaces that
are locally conformally flat, and besides the ellipsoid one knows
that the following types of hypersurfaces are conformally flat:
\begin{enumerate}
\item hypersurfaces of rotation;
\item centroaffine Tchebychev hypersurfaces with complete centroaffine metric 
and non-constant unimodular support function \cite{SIM-4}; 
\item decomposable hypersurfaces \cite{BI}.
\end{enumerate}
But one is far from a general  understanding of
the conformal properties in relative hypersurface theory, as there
are only  few results in special relative hypersurface theories.
Concerning conformal properties, this motivates a particular
interest in further investigations, and thus we consider
special relative hypersurfaces in the following subsections,
restricting to locally strongly convex relative hypersurfaces with
appropriate orientation such that the class $\mathfrak C$ is
(positive)  definite.

From section 8.2 recall the definition of the  trace free part $
\widetilde C$ 
of $C$. 
\begin{proposition}\label{prop-10.2}
On a non-degenerate hypersurface, let $Y$ be an arbitrary conormal field;
from $Y$ one can define the corresponding
metric $h$ and the projectively
flat connection $\nabla ^*$, and from this $C$, $T$ and finally
$\widetilde C$; we have:
\begin{enumerate}
\item $\widetilde C$ is  gauge invariant.
\item $\widetilde C = C(e)$.
\end{enumerate}\end{proposition}
%%%%%%%%%%%%%%%%%%%%%%%%%%%%%%%%%%%%%%
\subsubsection{Gauge invariant relative geometries}\label{sect-10.3.1}
See \cite{SIM-4}. We recall that 
the important relative hypersurface theories
are in fact gauge invariant; more precisely:

\begin{lemma}\label{lem-10.3}
One has that:
\begin{enumerate}
\item The centroaffine metric,  and thus its intrinsic geometry, 
is  gauge invariant.
\item The class  $\{c\cdot h$  $\arrowvert$   $h$  Blaschke metric, $0 \ne c \in 
\mathbb{R} \}$  is gauge invariant and thus also the intrinsic geometry of the
Blaschke metric  (modulo a non-zero constant factor).
\end{enumerate}\end{lemma}
%%%%%%%%%%%%%%%%%%%%%%%%%%%%%%%%%%%%%%%%%%%%%%%%%%%%%%%%%%%%%%%%%%
\subsection{Some special classes  of relative hypersurfaces}\label{sect-10.4}
%%%%%%%%%%%%%%%%%%%%%%%%%%%%%%%%%%%%%%%%%%%%%%%%%%%%%%%%%%%%%%%%%%
We list some special classes of hypersurfaces that are well known in relative 
hypersurface theory.
\subsubsection{Quadrics}\label{sect-10.4.1}
We have the following characterization of quadrics in terms 
of the gauge invariant cubic form \cite{SSV,SIM-4}.
\begin{theorem}\label{thm-10.4}
A relative hypersurface is a quadric if and only if 
$\widetilde C = 0$.
\end{theorem}
\subsubsection{Relative spheres}\label{sect-10.4.2}
Let $(x,Y,y)$ be a relative hypersurface;
it is called a {\it proper relative sphere} if, for some
$x_o \in \mathbb{R}^{n+1}$, we have $y=\lambda (x-x_o)$  for an appropriate
nowhere vanishing differentiable function $\lambda;$
 it is called an {\it improper relative sphere} if
$y$ is a constant transversal field. Furthermore, $(x,Y,y)$ is a relative sphere
if and only if $S^{\flat} = \lambda \cdot h = H \cdot h$, and the latter
relation implies $H=const.$  A relative sphere is proper if $H \ne 0$,
and it is improper if $H = 0$.

In the sense of the 
definition of relative spheres, any centroaffine hypersurface with centroaffine
normalization is a proper relative sphere, thus in the centroaffine
geometry the notion of  ``relative spheres'' is meaningless. In
Blaschke's geometry the relative spheres are called {\it affine spheres}.
For  proper affine spheres the Blaschke normalization and the
centroaffine  normalization coincide (modulo a non-zero constant factor),
thus the equation $T(c)=0$ characterizes proper affine spheres
within the centroaffine geometry.
The class of affine spheres is so large that one is far from any
local classification. Under additional assumptions there are
partial local and global classifications \cite{LI-S-Z-1,VRA}.
%%%%%%%%%%%%%%%%%%%%%%%%%%%%%%
%%%%%%%%%%%%%%%%%%%%%%%%%%%%%%
\subsubsection{Extremal hypersurfaces}\label{sect-10.4.3}
For any non-degenerate  hypersurface
$x$ which has a given conormal $Y$ the area functional
of its pseudo-Riemannian volume form $\omega (h)$, on a domain $D$
with compact support, is given by
$$ \mathfrak{B} := \int_D \omega (h). $$
In case a hypersurface is a critical point of the functional it
satisfies the Euler-Lagrange equation;  then the    hypersurface is
called an {\it extremal  hypersurface}.
%%%%%%%%%%%%%%%%%%%%%%%%%%%%%%%%%%%%%%%%
\subsubsection{Equiaffine extremal hypersurfaces}\label{sect-10.4.4}
We use the notation from section \ref{sect-8.3}.
The Euler-Lagrange equation takes the form
$nH(e):=Tr(S(e)) =0 $   in the  Blaschke geometry;  it is a PDE of fourth order.
The expression for the second variation of
the area functional is very complicated. E. Calabi \cite{CAL-2} proved:

\begin{theorem}\label{thm-10.5}
On locally strongly convex, extremal  hypersurfaces, any of the
following conditions (1) and (2) implies that the second
variation is negative; in this case the affine extremal hypersurfaces
are called {\it affine maximal}:
\begin{enumerate}
\item $n=2$.
\item $n \geq 2$ and $x$ can be represented as a graph.
\end{enumerate}\end{theorem}

\section{$W$-decomposition and relative hypersurfaces}\label{sect-11}
\subsection{$W$-decomposition and integrability conditions}
One of the integrability conditions in relative hypersurface theory
is given by the projective flatness of $\nabla^*$;  for $n\ge 3$,
Lemma 9.3 and Lemma 10.1 imply: 
 $$\pi_5(R^*) = \pi_6(R^*) =  \pi_7(R^*) = 0.$$ 
\subsection{Characterization of hyperquadrics}\label{sect-11.2}
Recall the definition of  $\widetilde C$ from section \ref{sect-8.2}
and the fact that
a hypersurface is locally strongly convex if and only if
the conformal class of relative metrics is (positive)  definite.

\begin{theorem}\label{thm-11.1}
 Consider a relative  hypersurface. Then:
\begin{enumerate} 
\item The expression 
$ \widetilde{\gamma}^i_{jkl} := \widetilde{C}^h_{jl} \widetilde{C}^i_{hk} - 
\widetilde{C}^h_{jk} \widetilde{C}^i_{hl}
$
defines an algebraic, gauge invariant curvature operator.
\item The Blaschke geometry  in the Codazzi structure has the following property: 
$$ \mathcal{R}(e) +  \mathcal{R}^*(e) - 2  \mathcal{R}(h(e)) = \widetilde{\gamma}\,.$$
\item In case the conformal class is (positive) definite,
we have the equivalences:
\begin{enumerate} 
\item  $ \widetilde{\gamma}$ vanishes identically.
\item  $ \widetilde{C} = 0$.
\item $ R(e) + R^*(e) = 2 R(h(e))$.
\item The gauge invariant connections in the Blaschke
geometry coincide: $$\nabla (e) = \nabla(h(e)) = \nabla^* (e) \,.$$
\item[(e)] $x$ is a hyperquadric.
\item[(f)] $\pi_1(R) = n(n-1)\,\pi_1(R(h)).$ 
\end{enumerate}
\end{enumerate}\end{theorem}

\begin{proof} We restrict to some remarks, as the proof is routine.
 In (3), the equivalences of (b), (d), (e) are true
for any relative  hypersurface.  In the Blaschke geometry, (f) yields if and only if the Pick invariant $J$ satisfies $n(n-1)J = n(n-1)(\kappa - H) =  n(n-1)\kappa   - \tau =0.$
As $g$ is positive definite, this is equivalent to 
$\tilde{C} = C = 0$ (see Theorem  \ref{thm-10.4}).
\end{proof}
%%%%%%%%%%%%%%%%%%%%%%%%%%%%%%%%%%%%%%%%%%%
\subsection{Characterization of relative and affine spheres}\label{sect-11.3}
We recall the following result from \cite{SSV}, Theorem 6.3.5.2.

\begin{theorem}\label{thm-11.2} Let $x: M \rightarrow \mathbb{R}^{n+1}$ be a 
centroaffine hypersurface of dimension $n>2$  with relative
normalization $(Y,y)$ and induced conjugate triple $(\nabla, h, \nabla^*)$. Then the 
relative normalization  coincides with the 
centroaffine normalization if and only if $\nabla$ is projectively flat.
\end{theorem}

This together with Section \ref{sect-11.5} above gives:

\begin{theorem}\label{thm-11.3}
Let $x: M \rightarrow \mathbb{R}^{n+1}$ be a 
centroaffine hypersurface of dimension $n>2$  with relative
normalization $(Y,y)$ and induced conjugate triple $(\nabla, h, \nabla^*)$. Then 
we have the 
equivalence of the following properties:
\begin{enumerate}
\item $\nabla$ is projectively flat.
\item $Ric = Ric^*$.
\item $\pi_2({R})=0$.
\item $\pi_2({R}^*) =0$.
\item $n\cdot Ric = \tau \cdot g$, that means  $(M,\nabla,h)$ is equiaffine Einstein.
\item $n\cdot Ric^* = \tau \cdot g$.
\item $\mathcal{R}  = \mathcal{R}^*$.
\item $-\,n(n-1)\cdot {R} = \tau \cdot (h \wedgeo  h)$.
\item The relative normalization  coincides with the 
centroaffine normalization.
\item $(x,Y,y)$ is a relative sphere.
\end{enumerate}\end{theorem}

If we apply Theorem 6.3.5.2 in \cite{SSV} and the foregoing Theorem we get:

\begin{theorem}\label{thm-11.4} The following assertions are
equivalent for a  relative  hypersurface $(x,Y,y)$ in dimension
$n \ge 3$.
\begin{enumerate} 
\item  $(x,Y,y)$ is a relative sphere.
\item $\nabla$ is projectively flat.
\item $B:= (n-1)Ric^* + Ric - \tau \cdot h  = 0$.
\item $\pi_2(R) = 0$.
\item  $\pi_2(R^*) = 0$.
\end{enumerate}\end{theorem}

\begin{remark}\label{rmk-11.5}
\rm The foregoing  Theorems  reflects the importance of 
the symmetric
bilinear form $B= \rho_{13}(\pi_5({R}))$; for  a relative  
 hypersurface
one  calculates that
$$B=n(n-2)(S^{\flat} - H\cdot h)\,.$$
In particular,   
one can characterize affine spheres by (3) above in the Blaschke geometry. Recall Section \ref{sect-10.4.2}
and the fact that 
this class of hypersurfaces is very large.
\end{remark}

\begin{theorem}\label{thm-11.6}
Let $(x,Y,y)$ be a non-degenerate relative hypersurface.
Then one has ${R} \in W_5$ if and only if the hypersurface is
an improper relative hypersphere.
\end{theorem}
\begin{proof}   ${R} \in W_5$
implies  that
$$0=Ric(\mathcal{R}^*)= n(n-1)S^{\flat},$$
thus   $(x,Y,y)$ is an improper relative hypersphere.
The converse is trivial. \end{proof}

The results in Section \ref{sect-11.5} specialize to affine spheres.
Recall that  for  non-degenerate hypersurfaces the conormal connection
$\nabla^*$ is always projectively flat, and the projective flatness
is equivalent to the relation
$${R}^* = \pi_1({R}^*) \oplus \pi_2({R}^*)\,.$$
Moreover, we have the equivalences:

\begin{observation}\label{obs-11.7}
For a Blaschke hypersurface adopt the notation established above. Then:
\begin{enumerate}
\item 
$\pi_2({R}) = 0$ if and only if $x$ is an affine sphere; this is 
equivalent to the total symmetry of $\nabla(h)\,C$.
\item ${R}(g) = \pi_1({R}(g))$ is equivalent to
the fact that the Blaschke metric has constant sectional curvature.
\end{enumerate}\end{observation}
%%%%%%%%%%%%%%%%%%%%%%%%%%%%%%%%%%%%%%
\subsection{Characterization of affine maximal hypersurfaces}\label{sect-11.4}
 Let $x$ be a locally strongly convex Blaschke hypersurface.
%%%%%%%%%%%%%%%%%%%%%%%%%%%%%%%%%%%%%%%%%%%%%%%%%%%%%%%%%%%%%%%%%
\begin{proposition}\label{prop-11.8}
The following properties  are equivalent:
\begin{enumerate}
\item $x$ is affine maximal.
\item $\pi_1({R}) = 0$.
\item $\pi_1({R}^*) = 0$.
\end{enumerate}\end{proposition}
%%%%%%%%%%%%%%%%%%%%%%%%%%%%%%%%%%%%%%%%%%%%%%%%%%%%%%%%%%%%
\subsection{Characterization of classes of Blaschke hypersurfaces in
  terms of PDEs and curvature tensors}\label{sect-11.5}
In the foregoing section we characterized some special classes
of hypersurfaces in terms of their equiaffine curvature tensors. On the
other hand it is well known that some of these classes can be also
locally characterized in terms of PDEs for a graph representation.
We combine characterizations
in terms of decomposition results from section \ref{sect-10}
with known characterizations in terms of PDEs from \cite{LI-S-Z-1}.

Let $x: M \rightarrow \mathbb{R}^{n+1}$  be a hypersurface
with a local representation by a strongly convex graph function
$f: \Omega \rightarrow \mathbb{R}$ with   $f=f(x^1,...,x^n),$ where $\Omega$ is a 
domain
in $\mathbb{R}^n$ s.t. the origin  $O \in \Omega \cap x(M)$  lies in the tangent plane 
$T_ox(M) = \mathbb{R}^n$.
\begin{theorem}\label{thm-11.9} Let $f$ be the above graph function. Then
\begin{enumerate} 
%%%%%%%%%%%%%%%%%%%%%%%%%%%%%%
\item $x$  is a proper affine sphere 
\begin{enumerate}
\item 
if and only if the Legendre transform function $u=
u(\xi^1,...,\xi^n)$ of 
$f$ satisfies the PDE 
$$det\,\left(\frac{\partial^2 u}{\partial\xi^i \partial\xi^j}\right) =
(H\,u)^{-(n+2)}\; ,$$
\item
if and only if $\pi_2({R}) = 0$.
\end{enumerate}
%%%%%%%%%%%%%%%%%%%%%%%%%%%%%%
\item $x$  is an improper affine sphere with constant affine normal vector
(0,...,0,1) 
\begin{enumerate}
\item 
 if and only if the graph function satisfies the Monge-Amp\`{e}re  PDE 
$$ det\,\left(\frac{\partial^2 f}{\partial x^i \partial x^j} \right) =
1,$$
\item 
 if and only if $\pi_5({R}) = 0$.
\end{enumerate}
%%%%%%%%%%%%%%%%%%%%%%%%%%%%%%
\item $x$  is affine maximal   
\begin{enumerate}
\item 
if and only if the graph function 
$f$ satisfies the PDE 
$$\Delta \left(\left[det\,(\frac{\partial^2 f}{\partial x^i \partial x^j} 
)\right]^{\frac{-1}{n+2}}\right) = 0,$$
\item 
 if and only if $\pi_1({R}) = 0$. 
\end{enumerate}
\end{enumerate}\end{theorem}
\medskip
 Lemma \ref{lem-7.6}(4) shows that, in a similar way,  affine spheres  in Blaschke's geometry
can be characterized in terms of PDEs
for the cubic form,  using results from \cite{BNS}.

The following Corollary states modifications of two  well
known global results, namely the Theorem of Blaschke and Deicke
and the affine Bernstein problem in the version of Calabi; 
see \cite{LI-S-Z-1}.
\begin{corollary}\label{cor-11.10}
Let $x: M \rightarrow \mathbb{R}^{n+1}$ be a locally
strongly convex Blaschke hyper-surface.
\begin{enumerate}
\item If $M$ is compact and if $\pi_2(R) = 0$ then $x$ is a hyperellipsoid.
 \item If $n=2$, $(M,h)$ is complete,   and   $\pi_1(R) = 0$, then $x$ is an
elliptic paraboloid.
\end{enumerate}
\end{corollary}

\section*{Acknowledgments} 
Research of 
 P. Gilkey partially supported by DFG
PI-158/5-5 and by Project MTM2006-01432 (Spain).
 Research of S. Nik\v cevi\'c partially supported by DAAD
 (Germany), TU Berlin, Dierks von Zweck Stiftung (Essen, Germany), and Project 144032 (Srbija). 
Research of 
 U. Simon partially supported by DFG PI-158/5-5.

\end{document}